\documentstyle[]{article}
\def\importpicture#1by#2(#3){
\vbox to #2 {
    \hrule width #1 height 0pt depth 0pt \vfill \special{picture #3}}
}

\def\scaledimportpicture#1by#2(#3scaled#4){{
\dimen0=#1  \dimen1=#2
\divide\dimen0 by 1000 \multiply\dimen0 by #4
\divide\dimen1 by 1000 \multiply\dimen1 by #4
\picture \dimen0 by \dimen1 (#3 scaled #4)}}
\def\dfigure#1by#2(#3scaled#4offset#5:#6)
    {\medskip
     \vglue 2mm minus 2mm
     $$
       \hbox{
         \hglue#5
         {\scaledimportpicture #1 by #2 (#3 scaled #4)}
       }
     $$
     \par\goodbreak
     \vglue 2mm minus 2mm
     \medskip}

\def\hugeslash{\mbox{ 
\unitlength=1pt
 \begin{picture}(0,0)
 \thinlines
 \put(0,-8){\line(1,3){8}}
\end{picture}
}}

\def\qmod#1#2{{\hbox{}^{\displaystyle{#1}}}\!\big/\!\hbox{}_{
\displaystyle{#2}}}

\def\bigqmod#1#2{{\hbox{}^{\displaystyle{#1}}} \!\hskip-4pt{\hugeslash}\!\hskip4pt\hbox{}_{
\displaystyle{#2}}}

\def\catqmod#1#2{{\hbox{}^{\displaystyle{#1}}}\!\big/\hskip-4pt{\big/}\!\hbox{}_{
\displaystyle{#2}}}

\newfam\msbfam

\font\tenmsb=msbm10
\font\eightmsb=msbm10 at 8pt
\font\sevenmsb=msbm10 at 7pt
\font\fivemsb=msbm10 at 5pt

\textfont\msbfam=\tenmsb
\scriptfont\msbfam=\sevenmsb
\scriptscriptfont\msbfam=\fivemsb

\def\Bbb{\fam\msbfam\tenmsb}

\def\A{{\Bbb A}}

\def\C{{\Bbb C}}

\def\N{{\Bbb N}}

\def\Q{{\Bbb Q}}
\def\R{{\Bbb R}}
\def\Z{{\Bbb Z}}

\def\union{\mathop{\bigcup}}
\def\qed {\hfill\vrule height6pt width6pt depth0pt \bigskip}

\def\map{\longrightarrow}
\def\textmap#1{\mathop{\vbox{\ialign{
                                  ##\crcr
      ${\scriptstyle\hfil\;\;#1\;\;\hfil}$\crcr
      \noalign{\kern-0.1pt\nointerlineskip}
      \rightarrowfill\crcr}}\;}}
\def\textra#1{\mathop{\vbox{\ialign{
                                  ##\crcr
      ${\scriptstyle\hfil\; #1\; \hfil}$\crcr
      \noalign{\kern-0.1pt\nointerlineskip}
      \rightarrowfill\crcr}}\;}}

\def\textlmap#1{\mathop{\vbox{\ialign{
                                  ##\crcr
      ${\scriptstyle\hfil\;\;#1\;\;\hfil}$\crcr
      \noalign{\kern-1pt\nointerlineskip}
      \leftarrowfill\crcr}}\;}}

\newfam\meuffam

\font\tenmeuf=eufm10
\font\sevenmeuf=eufm7
\font\fivemeuf=eufm5

\textfont\meuffam=\tenmeuf
\scriptfont\meuffam=\sevenmeuf
\scriptscriptfont\meuffam=\fivemeuf

\def\germ{\fam\meuffam\tenmeuf}

\def\cg{{\germ c}}

\def\g{{\germ g}}
\def\hg{{\germ h}}

\def\kg{{\germ k}}

\def\mg{{\germ m}}

\def\rg{{\germ r}}

\def\tg{{\germ t}}

\def\Eg{{\germ E}}
\def\Fg{{\germ F}}

\def\Kg{{\germ K}}
\def\Lg{{\germ L}}

\def\Pg{{\germ P}}

\begin{document}
\def\pro{{\rm pr}}
\def\tr{{\rm Tr}}
\def\End{{\rm End}}
\def\Aut{{\rm Aut}}
\def\Spin{{\rm Spin}}
\def\U{{\rm U}}
\def\SU{{\rm SU}}
\def\SO{{\rm SO}}
\def\PU{{\rm PU}}
\def\GL{{\rm GL}}
\def\spin{{\rm spin}}
\def\u{{\rm u}}
\def\su{{\rm su}}
\def\so{{\rm so}}
\def\ub{\underbar}
\def\proj{{\rm pr}}
\def\pu{{\rm pu}}
\def\Pic{{\rm Pic}}
\def\Iso{{\rm Iso}}
\def\NS{{\rm NS}}
\def\deg{{\rm deg}}
\def\Hom{{\rm Hom}}
\def\Aut{{\rm Aut}}
\def\h{{\germ h}}
\def\Herm{{\rm Herm}}
\def\Vol{{\rm Vol}}
\def\pf{{\bf Proof: }}
\def\id{{\rm id}}
\def\i{{\germ i}}
\def\im{{\rm im}}
\def\rk{{\rm rk}}
\def\ad{{\rm ad}}
\def\h{{\bf H}}
\def\coker{{\rm coker}}
\def\dv{\bar\partial}
\def\Ad{{\rm Ad}}
\def\Lie{{\rm Lie}}
\def\RSU{\R SU}
\def\ad{{\rm ad}}
\def\dva{\bar\partial_A}
\def\da{\partial_A}
\def\p{\partial\bar\partial}
\def\sp{\Sigma^{+}}
\def\sm{\Sigma^{-}}
\def\spm{\Sigma^{\pm}}
\def\smp{\Sigma^{\mp}}
\def\oo{{\scriptstyle{\cal O}}}
\def\ooo{{\scriptscriptstyle{\cal O}}}
\def\sw{Seiberg-Witten }
\def\pa{\partial_A\bar\partial_A}
\def\Dr{{\raisebox{0.15ex}{$\not$}}{\hskip -1pt {D}}}
\def\gr{{\scriptscriptstyle|}\hskip -4pt{\g}}
\def\subsetint{{\  {\subset}\hskip -2.45mm{\raisebox{.28ex}
{$\scriptscriptstyle\subset$}}\ }}
\def\nr{\parallel}
\def\ra{\rightarrow}
\newtheorem{sz}{Satz}[section]
\newtheorem{thry}[sz]{Theorem}
\newtheorem{pr}[sz]{Proposition}
\newtheorem{re}[sz]{Remark}
\newtheorem{co}[sz]{Corollary}
\newtheorem{dt}[sz]{Definition}
\newtheorem{lm}[sz]{Lemma}

\title{Comparing virtual fundamental classes: Gauge theoretical   Gromov-Witten
invariants for toric varieties}
\author{Ch. Okonek$^*$  
\and A. Teleman\thanks{Partially supported by: EAGER -- European Algebraic
Geometry Research Training Network, contract No HPRN-CT-2000-00099 (BBW
99.0030), and by SNF, nr. 2000-055290.98/1}}
\maketitle
\setcounter{section}{-1}

\section{Introduction}

Perhaps one of the most important mathematical results obtained at the end of the
last century is a principle which states that the  algebraic geometric notion of
{\it stability} is closely related to the global analytical notion {\it
solution of a Hermite-Einstein type differential equation}.  The first
form of this   principle, usually called the Kobayashi-Hitchin correspondence, 
was  proved and applied in a spectacular way by S. Donaldson [D1], [D2],
[DK] in the case of vector bundles; later it was generalized to a large class
of similar situations (see for instance [Hi], [Bra], [Mu1], [LT]). In
general, a   Kobayashi-Hitchin correspondence establishes an isomorphism
between a    moduli space of   stable algebraic geometric objects and a moduli
space of solutions of a certain (generalized) Hermite-Einstein equation.

Another fundamental concept introduced at the end of the last century was the
notion of {\it virtual fundamental class}. Roughly speaking, this theory allows
to endow {\it oversized} moduli spaces with a homology class with closed supports
(or a class in its Chow ring) whose degree equals the expected dimension of the
moduli space {\it in a canonical way}.

In the algebraic geometric framework this notion was introduced by
Behrend-Fantechi [BF] generalizing ideas from Fulton [F].  Another version of
this concept is due to Li-Tian [LiT]. 

There is also an analogous concept in the differential geometric framework: the
rigorous formalism   developed by Brussee [Br] allows to endow every  moduli
space associated to a  gauge theoretical problem of {\it Fredholm type} with a
canonical Cech homology class with closed supports, whose degree equals the
expected dimension of the moduli space.

We believe that:\\ \\
{\it The moduli space associated to a    gauge
theoretical problem of Fredholm type admits a canonical Behrend-Fantechi virtual
class as soon as  all   data are algebraic. Moreover, in this case the
Kobayashi-Hitchin correspondence maps the gauge theoretical virtual
fundamental class    to its canonical Behrend-Fantechi virtual fundamental
class.}\\

The importance of such a statement is obvious: it gives a universal comparison
principle not only for moduli spaces, but also for a large class of {\it
invariants} defined within the two categories.

The concept of "gauge theoretical problem of Fredholm type" is {\it much more
general} than one would think. For instance, in [OT2] we showed that the
vortex problem for line bundles on complex surfaces \ub{is} of Fredholm type,
although the elliptic deformation complex at a solution has a non-trivial 
degree 3 term. 

Note that a rigorous proof of our conjecture cannot be easy: whereas the
definition in gauge theory uses Sobolev completions, Fredholm operators and 
Cech homology with closed supports, the definition in the algebraic category uses
sheaves in the \'etale topology, derived categories, Deligne-Mumford
stacks, cotangent complexes and Chow rings!

The purpose of this paper is neither to   give a precise general form of our
conjecture, nor to speculate on a possible proof strategy.  We will come back to
the general case in future works.

Our purpose here is to illustrate this principle in an interesting  concrete
situation for which we will make a precise statement and give a complete
proof.\\

Our concrete situation is the moduli problem considered in the definition of
the {\it twisted gauge theoretical Gromov-Witten invariants associated with the
symplectic factorization problem used to construct   complete toric
varieties.}\\

 In our previous paper [OT2] we introduced a new type of gauge
theoretical  Gromov-Witten invariants which    generalize the so
called "Hamiltonian Gromov-Witten invariants" introduced independently
in [Mu] and [CGS] (see also [CGMS]).

Our invariants are associated  with triples $(F,\alpha,K)$, where $(F,\omega,J)$    
is an almost K\"ahler manifold, $\alpha$ is a $J$-holomorphic action of a compact
Lie group $\hat K$ on $F$, and $K$ is a closed normal subgroup of $\hat K$ which
leaves the symplectic form $\omega$ invariant.

In the quoted article we
called such triples {\it symplectic factorization problems with
additional symmetry}, because  we only consider   symplectic
quotients  with respect to the normal subgroup $K$, whereas the
manifold $F$ was endowed with the action of a larger Lie group $\hat K$,
which will act (in general)   non-trivially on the symplectic
$K$-quotients of $F$.

The quotient group    $K_0:=\hat
K/K$ (called the parameter symmetry group or the {\it twisting} group) plays an
important role in our approach [OT2]. The formalism developed in  [CGS], [CGMS]
corresponds to the  case when $K_0$ is trivial.

Let us  denote by $\pi:\hat K\ra K_0$   the canonical projection.
Our invariants are obtained by evaluating canonical cohomology classes
on the virtual fundamental class of the moduli space of solutions of a
certain       generalized vortex equation, which depends on the choice
of:\vspace{2mm}
\\ 1.   a system of (discrete) {\it topological 
parameters}, namely a triple $(Y,P_0,\cg)$ consisting of a closed
oriented real surface
$Y$, a $K_0$-bundle $P_0$ on $Y$,   an equivalence class
$\cg$ of pairs $(\hat P\textra{\lambda} P_0,\hat h)$ formed by a morphism of
type $\pi$, and a homotopy class $\hat h$ of sections in the associated bundle
$E:=\hat P\times_{\hat K} F$.\vspace{2mm}\\
2.  a system of {\it continuous parameters}, namely a triple
$(\mu,g,A^0)$ consisting of a Riemannian metric $g$ on $Y$, a
connection $A^0$ in $P_0$, and a
$\hat K$-equivariant moment map $\mu$ for the $K$-action of $F$.\\

The  first purpose of this  paper is to apply this general set up  in
order to introduce    the resulting invariants  rigorously in the
following  important special case: $F=\C^r$, $\alpha$ is the natural
action $\alpha_{\rm can}$ of
$\hat K=[S^1]^r$       on
$F$, and $K$ is the kernel $K_w$ of an epimorphism $w: [S^1]^r\ra 
K_0=[S^1]^m$, hence a compact (but possibly non-connected) abelian group of
dimension $r-m$.  Therefore, the invariants we introduce and study should be
  called {\it $K_w$-equivariant, $[S^1]^m$-twisted gauge
theoretical Gromov-Witten invariants of the affine space $\C^r$.}

If the adiabatic limit conjecture is true (see [G], [CGS], [GS]), these
invariants should be related to the {\it twisted} Gromov-Witten
invariants of   toric varieties. These twisted Gromov-Witten
invariants, which were introduced in [OT2], are natural
generalizations of the Gromov-Witten invariants in the sense of
Ruan [R]. They are obtained by replacing the moduli spaces of (pseudo)holomorphic
morphisms $Y\ra F$ in   Ruan's definition of  
Gromov-Witten invariants  by   moduli spaces of (pseudo)holomorphic sections in
$F$-bundles over $Y$ with a fixed structure group $K_0$. In other words, we
replace the gauged linear  
sigma models introduced by Witten   [W3], and further investigated by
Morrison und Plesser [MP], by {\it
$K_0$-twisted gauged linear  
sigma models} [OT2].\\

Our first result  is a Kobayashi-Hitchin correspondence which also
gives an explicit  complex geometric interpretation for the virtual
fundamental class   of  the moduli space.  More precisely:

Let $V\in M_{m,r}(\Z)$ be an integer matrix of
rank
$m$,   let
$w:[S^1]^r\ra [S^1]^m$ be the corresponding epimorphism, and  put
$K_w:=\ker(w)$.   Suppose
that the columns of
$V$ are primitive and that 
$$^tV(\R_m)\cap\{(t_1,\dots,t_r)\in\R_r|\ t_i\geq 0\}=\{0\}\
.\footnote{This condition assures that the symplectic quotients of
 {\eightmsb C}$^r$ by $K_w$ are compact.}
$$

Let $t\in \kg_w^{\vee}=i\coker(^tV\otimes\id_\R)$ be a regular value of
the standard moment map
$\mu_w$ of the
$K_w$-action on $\C^r$.  Under these assumptions, the K\"ahler quotients
$$X_\tg:=\catqmod{\C^r}{_{[\mu_w+i\tg]}K_w}
$$
is a compact toric variety with a natural orbifold structure.

Let $Y$ be a closed oriented surface. We fix a $[S^1]^r$-bundle $\hat P$,
a
$[S^1]^m$-bundle
$P_0$ and a $w$-morphism   $\hat P\textmap{\lambda}P_0$ over $Y$.
Let   $g$ be a Riemannian metric on $Y$, and let $A^0$ be
a connection on $P_0$. Then
\begin{thry}\hfill{\break}\vspace{-6mm}
\begin{enumerate}
\item (complex geometric interpretation) The moduli space ${\cal M}$ of
  solutions of the vortex equation associated with the  data
$(\lambda,(\mu_w+i\tg,g,A^0))$ is a toric fibration over an abelian 
variety $P$ of dimension $g(Y)(r-m)$.
\item (embedding theorem) The moduli space ${\cal M}$ can be identified
with  the vanishing locus of a   section $\sigma$ in an explicit 
holomorphic   bundle  
${\cal E}$ over the total space of a   locally trivial holomorphic toric fibre
bundle
$T$ over $P$.
\end{enumerate}
\end{thry}

The standard fibre $\Phi$ of the toric fibre  bundle $T$ is a toric
orbifold which can be obtained as a K\"ahler quotient of a suitable
complex vector space by the same group $K_w$ used to construct the
toric variety $X$. In the case $g(Y)=0$ this has already been observed in
[MP]. This fibre
$\Phi$ is a smooth manifold if the  K\"ahler   quotient $X$  was
smooth. In this case one gets an identification $\iota:{\cal M}\ra
Z(\sigma)$ between 
${\cal M}$ and the subspace cut out by a holomorphic section $\sigma$ in
a holomorphic vector bundle
${\cal E}$ over a smooth complex manifold $T$.

This allows us to endow ${\cal M}$ with a distinguished homology class
of degree $\dim(T)-\rk({\cal E})$, namely the {\it algebraic geometric}
virtual fundamental class  of the triple $(T,{\cal E},\sigma)$ in the
sense of Fulton [Fu]. In the general case,
$T$ can be regarded as a Deligne-Mumford stack, so one can still endow
${\cal M}$ with a rational virtual fundamental class in the sense of
Behrend-Fantechi [BF].\\

{\it Unfortunately, the construction of $T$, ${\cal E}$, $\sigma$ and of
the embedding ${\cal M}\hookrightarrow T$ is \ub{not} \ub{canonical}: it
depends on the choice of a system of  sufficiently ample divisors
$(D_j)_{1\leq j\leq r}$ on
$Y$, and it is not  clear whether the virtual fundamental
class obtained in this way is independent of this choice.}
\\

In order to prove this, one
should check that different choices of systems of ample divisors lead to
the same {\it perfect obstruction theory} in the sense of
Behrend-Fantechi.

Our next result states:
\begin{thry} (comparison theorem)  The algebraic geometric virtual
fundamental class of
${\cal M}$ induced by the  identification $\iota$ coincides with its
gauge theoretical virtual fundamental class.
\end{thry}

Recall that the gauge theoretical virtual fundamental class is obtained
by regarding
${\cal M}$ as the vanishing locus of a Fredholm section in a Banach
bundle over  a Banach manifold. The precise formalism was developed in 
[Br] (see also [OT2]). The obtained virtual class \ub{\it is} \ub{\it
canonical}, i. e. it does not depend on any additional   choices
 besides the parameters involved in the construction of the moduli
space.

The above theorem is very important, because it allows the
computation of the gauge theoretical (Hamiltonian) Gromov-Witten
invariants of   toric manifolds with purely algebraic geometric methods.
Since the pair
$(T,{\cal E})$ comes with a natural $[\C^*]^m$-action, and the section
$\sigma$ is equivariant, one can apply the localization theorem of   
Graber-Pandharipande [GP] for the Behrend-Fantechi virtual fundamental
class. Explicit computations can be found in [Ha1], [Ha2].

In this way, one can compare -- in the case of toric manifolds -- the
gauge theoretical (Hamiltonian) Gromov-Witten invariants with the
standard (Kontsevich)  Gromov-Witten invariants and check the adiabatic
limit conjecture for this class of manifolds. Explicit computations of
standard Gromov-Witten invariants for toric varieties can be found in
[Sp].

And now a word about the proof of the comparison theorem, Theorem 0.2
(see section 5):

  The configuration space of our gauge
theoretical problem is the  product of two factors: a space of
connections (or semiconnections) and a space  of sections. Our method is
based on the following new idea: we complete the space of sections in the
configuration space with respect to a {\it very weak Sobolev norm}, such
that   meromorphic sections with first order  poles in finitely many
simple points  become elements in our Sobolev completion. The spaces of
(semi)connections and the gauge group are
  completed as usual  with respect to    $L^2_k$ Sobolev norms.  This
asymmetric Sobolev completion  of the  configuration space   allows us
to pass from the gauge theoretical to the algebraic geometric framework. Of
course, one has to check that the new completed configuration space
leads to the same gauge theoretical virtual fundamental class as the
standard completion. The comparison Theorem 0.2 will  follow from
Brussee's {\it associativity principle} for virtual fundamental classes
associated with Fredholm sections [OT2].

We believe that this method (weakening the
Sobolev norm on the spaces of sections)  can be adapted to a very  large
class of similar problems. It can be used to show that many
Kobayashi-Hitchin correspondences (which relate gauge theoretical  to
algebraic geometric moduli spaces) map  the gauge theoretical virtual
fundamental class onto the algebraic geometric one. 

Another interesting application of this technique will be considered in
[DOT].

\section{Toric varieties as symplectic quotients}

Any epimorphism $w: [S^1]^r\map 
[S^1]^m$ is determined by the associated Lie algebra morphism, hence by a  linear
map $v:\R^r\ra
\R^m$ given by an integer matrix   $V=(v^i_j)_{\matrix{\scriptstyle1\leq i\leq m 
\vspace{-1.5mm}\cr\scriptstyle1\leq j\leq r }}\in M_{m,r}(\Z)$  of rank
$m$; one has
$$w (e^{ i t^1},\dots,e^{ i t^r})=(e^{i v^1_j t^j},\dots,e^{iv^m_j t^j})\ .
$$

We are not interested in {\it all} epimorphisms $[S^1]^r\ra 
[S^1]^m$ as above, but only in those epimorphisms $w$ with the property that  
the symplectic quotients of the form $\C^r/_\mu K_w$ are   compact, because  in
this case, the corresponding invariants   should be related to the 
{\it twisted   Gromov-Witten invariants of   toric varieties } 
[OT2].\\

Therefore we shall assume that  $w$ verifies the following properties:\\ \\
{\bf P$_1$: } {\it For every $j\in\{1,\dots,r\}$, the column  $v_j\in
\Z^m$ is primitive, i. e. it is a generator of the semigroup $\Z^m\cap
\R_{\geq 0} v_j$.}
\\ \\
{\bf P$_2$: } {\it  $\R_r^{\geq 0}\cap \im(v^*)=\{0\}$ in  the dual 
space  
$\R_r$ of $\R^r$. }
\vspace{2mm}\\
Here we used the notation 
$\R_r^{\geq 0}:=\{(t_1,\dots,t_r)\in\R_r|\ t_i\geq 0\}$.\\

Note that the second property is equivalent to
$$
\sum_{j=1}^r \R_{\geq 0}v_j= \R^m\ .
$$
Applying the functor $\Hom(\ \cdot\ , S^1)$ to the exact sequence
$$0\map K_w\map [S^1]^r\map [S^1]^m\map 0
$$
we get the exact sequence
$$0\map \Z_m\textmap{V^*} \
 \Z_r\map\Hom(K_w,S^1)\map 0\ .
$$
This shows that one has   natural identifications 
$$K_w=\Hom(\coker(V^*),S^1)\ ,\ \kg_w =i\ker(v)\
,\ \kg_w^\vee=i \coker(v^*)$$

The standard moment map $\hat\mu:\C^r\ra i  \R_r=\hat \kg^{\vee}$ of
the canonical action of
$[S^1]^r$ on
$\C^r$ is given, with respect to the  standard dual basis, by  
$$\hat \mu(z^1,\dots,z^r)=-\frac{i}{2}(| z^1|^2,\dots,|z_r|^2)\ .
$$
  The standard moment map $\mu_w:\C^r\ra  \kg^\vee_w=i \coker(v^*)$ of
the
$K_w$-action on
$\C^r$ is defined by
$$\mu_w(z^1,\dots,z^r)=-\frac{i}{2}p_v(| z^1|^2,\dots,|z_r|^2)\ ,
$$
where $p_v$ is the canonical projection $\R_r\ra \coker(v^*)$.
The image of $\mu_w$ is the convex set
$$A^+_v=-i p_v (\R^{\geq 0}_r)\subset i \coker(v^*)\ .
$$
An immediate consequence of the assumption ${\bf P}_2$ is
\begin{lm}  \hfill{\break}\\
i)  There is a  constant $c>0$ such that
$$\nr \mu_w(z)\nr^2 \geq  c \nr  z \nr^4\ .
$$
ii) All symplectic quotients
$$\qmod{\mu_w^{-1} (-i\tg)}{ K_w}\ , \ \ -i\tg \in A^+_v
$$
are  compact.
\end{lm}

The symplectic quotients which correspond to regular values of  the
moment map $\mu_w$ are projective toric varieties  with (at most) orbifold
singularities.
 
Conversely, let  $J\subset\{1,\dots,r\}$, and let $\Sigma$ be a {\it
complete,  simplicial} fan of strictly convex rational polyhedral cones in
$\R^m$ whose 1-skeleton  $\Sigma(1)$ is
$$\Sigma(1)=\{\R_{\geq 0} v_j |\ j\in J\}\ .$$

Let $a=(a_1,\dots,a_r)\in \R_r$.  For every
$\sigma\in\Sigma$ we define the functional  $f_\sigma^a\in\ 
\langle\sigma\rangle^\vee$ by requiring
$$\langle f_\sigma^a, v_j\rangle=-a_j\ \ {\rm if}\ \R_{\geq 0} v_j\ {\rm
is\ a \ face\ of\ }\sigma\ .
$$
The system $(f_\sigma^a)_{\sigma\in \Sigma}$ depends only on $(a_j)_{j\in
J}$ and it defines a continuous piecewise linear function $f^a$ on
the support
$|\Sigma|=\R^m$ of
$\Sigma$.  
We put:
$$K(\Sigma):=\{p_v(a)\ \mid a_i\geq 0, \ \  \langle f^a_\sigma,
v_j\rangle\geq -a_j\ \forall \sigma\in\Sigma\ ,\  \forall j\in
\{1,\dots,r\}\}
$$
$$K_0(\Sigma):=\{p_v(a)\in  K(\Sigma)\ \mid \langle f^a_\sigma,
v_j\rangle> -a_j\ \forall \sigma\in\Sigma\ ,\  \forall j\in
\{1,\dots,r\} $$
$$\ {\rm for\ which}\ \R_{\geq	0} v_j\ {\rm is\ not\ a\ face\
of\ }\sigma\}
$$

The conditions in the two definitions  depend 
only   on the class
$[a]$ modulo
$\im(v^*)$, because changing $a$ by an element of the form $v^*(f)$,
$f\in\  \R_m$    modifies all  maps $f_\sigma^a$ 
 by the same linear functional $f$.

When $p_v(a)\in K(\Sigma)$ ($K_0(\Sigma)$) the
piecewise linear map $f^a$ is convex (strictly convex) on
$|\Sigma|$ (see [Co]).  We denote by $k(\Sigma)$
(respectively $k_0(\Sigma)$) the cone of (strictly) convex $\Sigma$-linear maps
on
$|\Sigma|$.  One has obvious surjective maps $K(\Sigma)\ra k(\Sigma)$,
$K_0(\Sigma)\ra k_0(\Sigma)$.  Note that the cones $k(\Sigma)$, $k_0(\Sigma)$
depend only on the simplicial fan $\Sigma$, whereas $K(\Sigma)$,
$K_0(\Sigma)$ also depend on the rays $\R_{\geq 0 }v_j$ which are not faces
of $\Sigma$.

Recall that  every complete simplicial fan $\Sigma$ in $\R^m$ with
$$\Sigma(1)\subset\{\R_{\geq	0}v_1,\dots,\R_{\geq	0}v_r\}$$
defines an associated     compact
toric variety $X_\Sigma$  as follows: Set
$$U(\Sigma)=\{z\in\C^r|\  \exists\sigma\in\Sigma \ {\rm such\ that}\
z^j\ne0\ 
\forall   j\in\{1,\dots,r\}\ {\rm for\ which}  $$
$$\ \R_{\geq 0}v_j \ {\rm
is\ not\ a\ face\ of}\ \sigma  
\}\ .
$$
Then there is a geometric quotient
$$X_\Sigma:=\catqmod{U(\Sigma)}{K_w^\C}
$$
 and   this quotient is a compact algebraic (not necessarily projective) 
variety with a natural orbifold structure. The variety
$X_\Sigma$ is projective if and only if $k_0(\Sigma)$ (or equivalently
$K_0(\Sigma)$) is non-empty. In this case  one has a canonical
epimorphism $\coker(v^*)\ra H^2(X_\Sigma,\R)$
and, under this epimorphism, $K_0(\Sigma)$ 
is mapped onto the K\"ahler cone of the
orbifold $X_\Sigma$, which can be identified with $k_0(\Sigma)$. 
  We refer to [Co], [Gi2] for more details and the following theorem
\begin{thry} Let $\Sigma$ be a complete simplicial fan $\Sigma$ in $\R^m$
with
$$\Sigma(1)\subset\{\R_{\geq	0}v^1,\dots,\R_{\geq	0}v^r\}\ .$$ 
i)  For every $\tg\in K_0(\Sigma)$, the set of semistable points with
respect to the moment map $\mu_w+i\tg$ coincides with the correponding
set of stable points, and the symplectic quotient
$\qmod{\mu^{-1}_w(-i\tg)}{K_w}$ can be identified as a complex orbifold
with the projective toric variety $X_\Sigma$.\\
ii)  (The GKZ decomposition) The cone $p_v (\R^{\geq 0}_r)$ can be
decomposed as a union of subcones whose interiors are pairwise disjoint:
$$p_v (\R^{\geq 0}_r)=\union\limits_{\matrix{\scriptstyle\Sigma\ {\rm 
complete\ simplicial\ fan\ in\
}\R^m\vspace{-1.5mm}\cr\scriptstyle\Sigma(1)\subset\{\R_{\geq
0}v_1,\dots, \R_{\geq 0}v_r\}\vspace{-1.5mm}\cr
\scriptstyle k_0(\Sigma)\ne\emptyset}}  K(\Sigma)\ ,\ 
$$
iii) The nonempty open subcones $-iK_0(\Sigma)$ are  the
connected components of the complement  of the   critical locus ${\rm
Crit}(\mu_w)$ of $\mu_w$ in $\im(\mu_w)=-ip_v (\R^{\geq 0}_r)$.  Moreover, one
has
$${\rm Crit}(\mu_w)=\union_{\matrix{\scriptstyle
J\subset\{1,\dots,r\}\vspace{-1.5mm}\cr
\scriptstyle | J|=m+1}} 
\mu_w(Z_I)\ ,
$$
where $Z_J\subset \C^r$ is the subspace defined by the equations $z^j=0$,
$j\in J$.
\end{thry}

\section{Moduli spaces of toric vortices and the associated invariants}

\subsection{Moduli spaces of toric vortices}
 
Let  $Y$ be a closed connected oriented   real surface. The data
of a
$[S^1]^m$-bundle $P_0$ on $Y$ is equivalent to the data of a system 
$(L_i^0)_i$ of $m$ Hermitian line bundles
$L_i^0$ on $Y$.  

Fix a rank $m$ integer  matrix $V\in M_{m,r}(\Z)$    with primitive columns and
let $w:[S^1]^r\ra [S^1]^m$ be the associated epimorphism. The data of a
$w$-morphism
$\lambda:\hat P\ra P_0$ of principal bundles is equivalent to the data of
a system
$(L_j)_j$ of $r$ Hermitian line bundles and a   system  of
$m$ unitary isomorphisms
$$\lambda_i:\otimes_{j=1}^r [L_j^{\otimes v^i_j}]\map L_i^0\ .
$$

Following the general strategy explained  in [OT2], we consider the following
moduli problem: 
\\

Fix a Riemannian metric $g$ on $Y$, a 
system of Hermitian connections $A^0=(A_i^0)_{1\leq i\leq m}$  on
$(L_i^0)_{1\leq i\leq m}$, and an element
$\tg\in\
\coker (v^*)$. Classify all systems
$(A_j,\varphi_j)_{1\leq j\leq r}$, consisting of \\ \\
 i)  a connection  $A_j$   on
$L_j$  for every $j$, $1\leq j\leq r$ such that
$$\lambda_i(\otimes_{j=1}^r [A_j^{\otimes v^i_j}])=A_i^0\ ,\ \forall  \
i\in\{1,\dots,  m\} ,\eqno{(1)}$$
ii)  an $A_j$-holomorphic section $\varphi_j$   in
$L_j$ for every
$j\in\{1,\dots, r\}$\ , \\  \\
such that $(A_j,\varphi_j)_{1\leq j\leq r}$ solves the vortex-type
equation
$$p_v\left[(i\Lambda_gF_{A_j}-2\pi\deg(L_j)+\frac{1}{2}|\varphi_j|^2)_j\right]=
\tg\ .\eqno{(V^\tg_g)}
$$

Two such systems are considered equivalent if they are in the same orbit with
respect to  the natural action of the gauge group ${\cal G}:={\cal
C}^{\infty}(Y,K_w)$, i. e.
$${\cal G}=\{(f_1,\dots f_r)\in{\cal C}^\infty(Y,S^1)^r\ |\ \prod\limits_{j=1}^r 
f_j^{v^i_j}=1\ \forall i\in\{1,\dots,m\}\} \ .
$$
\begin{lm}\hfill{\break}
i) If the equation $(V_g^\tg)$ has a solution $(A_j,\varphi_j)_j$, then
$$
\tg   \in\union\limits_{\matrix{\scriptstyle\Sigma\
{\rm  complete\ simplicial\ fan\ in\
}\R^m\vspace{-1.5mm}\cr\scriptstyle\Sigma(1)\subset\{\R_{\geq
0}v^1,\dots, \R_{\geq
0}v^r\}\vspace{-1.5mm}\cr
\scriptstyle k_0(\Sigma)\ne\emptyset}}  K(\Sigma) \ .
$$
ii) If $\tg\in K_0(\Sigma)$ and  $(A_j,\varphi_j)_j$ is a solution  of
$(V^\tg_g)$, then   
$$(\nr \varphi_j\nr_{L^2})_{1\leq j\leq r}\in U(\Sigma)\ .$$ 
\end{lm}
\pf
Indeed, integrating $(V^\tg_g)$ over $Y$, we get
$$\tg \ {{\rm Vol}_g(Y)} =p_v\left[(\frac{1}{2}\nr
\varphi_j\nr_{L^2}^2)_j\right]\in p_v (\R_{r}^{\geq 0})\ ,
$$
which proves the first statement.  For the second, the same argument gives
$$-\frac{i}{2}p_v((\nr \varphi_j\nr_{L^2}^2)_{j})=-i\tg {{\rm Vol}_g(Y)}\ ,
$$
which implies that the system $(\nr \varphi_j\nr_{L^2})_{1\leq j\leq r}$
is semistable with respect to  the moment map $\mu_w+i\tg$.  The result
follows now from Theorem 1.2
\qed

The configuration space for our moduli problem is the product 
$${\cal A}:=\left[{\prod_{j=1}^r} {\cal A}(L_j) 
\right]_{A^0}^{(\lambda,v)}
\times\bigoplus_{j=1}^r A^0(L_j)\ ,
$$
where the first factor denotes the affine subspace of ${\prod_{j=1}^r} {\cal
A}(L_j)$ consisting of all systems of Hermitian connections $A=(A_j)_j$
satisfying  the relation $(1)$. Let
${\cal A}^{(H,V_g^\tg)}$ be the subspace of
${\cal A}$ cut out by  the integrabilty condition
$$\bar\partial_{A_j}\varphi_j=0\  ,\ j=1,\dots,r \eqno{(H)}
$$
and the vortex equation $(V_g^\tg)$. Our moduli space is the quotient
$${\cal M}_{(\tg,g,A^0)}(\lambda)=\qmod{{\cal A}^{(H,V_g^\tg)}}{\cal
G}\  ,
$$
and will be called {\it the moduli space of toric vortices} associated to
the data $(\tg,g,A^0,\lambda)$. 
\begin{pr}  The moduli spaces ${\cal M}_{(\tg,g,A^0)}(\lambda)$ are
compact.
\end{pr}
\pf  The proof uses the same argument as the demonstration of the compactness
Theorem 2.12 in [OT2].  The crucial point   is the
properness of the moment map which, in the present case, is stated in Lemma
1.1, $i)$.
\qed

Let ${\cal A}^*\subset{\cal A}$ (${\cal A}^{**}\subset{\cal
A}$) be  the space of points  with finite (trivial)
stabilizers  in the configuration space, and put
$${\cal B}^{**}:={\cal A}^{**}/{\cal G}\ ,\ {\cal B}^{*}:={\cal
A}^{*}/{\cal G}\ .$$
After suitable Sobolev completions ${\cal B}^{**}$ (${\cal B}^*$) 
becomes a Banach manifold (Banach orbifold). In both cases, the
local models are obtained in the usual way, by constructing local slices for  the
${\cal G}$-action in the configuration space (see for instance [OT1]).
\\

The maps $m_\tg:{\cal A}\ra
A^0(Y)^{\oplus m}$, $h :{\cal A}\ra \oplus_j A^{01}(L_j)$ given by
$$m_\tg(A,\varphi)=p_v\left[(i\Lambda_gF_{A_j}-2\pi\deg(L_j)+
\frac{1}{2}|\varphi_j|^2)_j\right] -
\tg \ ,\ h(A,\varphi)=(\bar\partial_{A_j}\varphi_j)_j
$$
are ${\cal G}$ - equivariant, hence they descend to real analytic
sections
$\mg_\tg$, $\hg$ in Banach bundles (respectively Banach orbifold 
bundles) over
${\cal B}^{**}$ (respectively ${\cal B}^*$). Moreover, the section
$(\mg_\tg,\hg)$ in the product bundle is Fredholm.

\begin{pr}  
If $K_w$ acts freely on $U(\Sigma)$ then the
toric variety $X_\Sigma$ is smooth, and  
the  moduli space
${\cal M}_{(\tg,g,A^0)}(\lambda)$ is a real analytic subspace of the
Banach manifold ${\cal B}^{**}$ for every $\tg\in K_0(\Sigma)$. It can be
identified with the vanishing locus of the Fredholm section
$(\mg_\tg,\hg)$.

In the general case,   the $K_w$-action on   $U(\Sigma)$ has
finite stabilizers and the moduli space ${\cal M}_{(\tg,g,A^0)}(\lambda)$ is a   
real analytic  suborbifold   of the infinite dimensional orbifold  ${\cal
B}^{*}$. It can be identified with the    real
analytic  suborbifold    cut out by the section
$(\mg_\tg,\hg)$.
\end{pr}
\pf  The point is that, if $\tg\in K_0(\Sigma)$, the stabilizer of a
solution $(A_j,\varphi_j)_j$ of $(V^\tg_g)$ with respect to the action
of the gauge group ${\cal G}$ coincides with the stabilizer of the point
$(\nr \varphi_j\nr_{L^2})_{1\leq j\leq r}\in U(\Sigma)$ with respect to
the Lie group $K_w$ (see Lemma 2.1).

Indeed, if a gauge transformation $f$ leaves the connection
system $(A_j)_j$ invariant, it must be a constant gauge transformation,
i. e. an element of $K_w$. On the other hand, a constant gauge
transformation   leaves  the section system
$(\varphi_j)_j$ invariant if and only if it leaves    the vector $(\nr
\varphi_j\nr_{L^2})_{1\leq j\leq r}\in U(\Sigma)$ invariant.

This shows that the gauge group acts freely (with finite stabilizers)
on the space of solutions of $(V^\tg_g)$ if $K_w$ acts freely (with
finite stabilizers) on $U(\Sigma)$. Therefore ${\cal
M}_{(\tg,g,A^0)}(\lambda)\subset {\cal B}^{**}$ in the  first case and   
${\cal M}_{(\tg,g,A^0)}(\lambda)\subset {\cal B}^{*}$ in the second.

In the second case, ${\cal M}_{(\tg,g,A^0)}(\lambda)$
has a natural orbifold structure whose local models are finite
dimensional real analytic spaces endowed with analytic actions of finite
groups. 
\qed

This proposition allows us -- at least in the case when the $K_w$-action on
$U(\Sigma)$ is free -- to endow the moduli space  ${\cal
M}_{(\tg,g,A^0)}(\lambda)$ with a virtual fundamental class in Brussee's sense
([Br], [OT2]):
$$[{\cal M}_{(\tg,g,A^0)}(\lambda)]^{\rm vir} \in H_{2e}({\cal
M}_{(\tg,g,A^0)}(\lambda),\Z)\ ,
$$
where $e=\sum_j\chi(L_j)-(r-m)\chi({\cal O}_Y)$ is the 
expected complex dimension of the moduli space. The obtained virtual
fundamental class {\it     does not depend on the chosen
Sobolev completions} (see the proof of  Proposition 5.4 and  
Remark 5.5). This is a rather non-trivial  fact which holds for a
large class of similar  gauge theoretical moduli spaces.

The definition of   virtual fundamental classes can be extended to the
orbifold case and then yields a  class in the rational
homology of the moduli space. This generalization will be treated in a
future work.

\subsection{Canonical cohomology classes on the moduli spaces of toric
vortices}

As always in gauge theory, one defines invariants by evaluating {\it
canonical cohomology classes} on the (virtual) fundamental class of the
moduli space. This general principle was applied in [OT2] in the case
of gauge theoretical Gromov-Witten invariants.  In the general case, the
canonical cohomology classes on the moduli space ${\cal M}$ associated
with a triple $(F,\alpha,K)$ and a system of parameters 
$((Y,P_0,\cg),(\mu,g,A^0))$ are products of classes
$\delta^\gamma(h)$ of the form
$$\delta^\gamma(h)=\Phi^*(\gamma)/h\ ,$$
where $\Phi$ is the universal section in the universal $F$-bundle over
the product ${\cal M}\times Y$, $\gamma\in H^*_{\hat K}(F,\Z)$ and $h\in
H_*(Y,\Z)$. The definition makes sense when the gauge group acts freely on
the space of solutions.  In the quoted article we showed that, in
general, these classes satisfy a set of tautological identities, so that
the relevant invariant can be regarded as a map on a quotient algebra
$\A$ of the graded algebra generated by the symbols
$\left(\matrix{\gamma\cr h}\right)$. The   algebra $\A$   depends only
on the homotopy type of the topological parameters $(Y,P_0,\cg)$.  
\\

In this section we describe  these  canonical cohomology 
classes  and the corresponding tautological relations explicitly in the special
case we are studying, namely in the case of the
 triple
$(\C^r,\alpha_{\rm can},K_w)$. Note first that the equivariant cohomology
ring
$H^*_{\hat K} (F,\Z)$
can be identified with the polynomial ring $\Z[c_1,\dots,c_r]$ in $r$ variables
of degree 2.

Consider an element $\tg\in\coker(v^*)$, and  suppose that
there exists a complete simplicial fan $\Sigma$ in $\R^m$ with $\Sigma(1)\subset\{\R_{\geq
0}v^1,\dots, \R_{\geq 0}v^r\}$ such that $\tg\in K_0(\Sigma)$.
\\ \\
\ub{Case} 1. The action of $K_w$ on $U(\Sigma)$ is free.\\

In this case the toric variety  
$X_\Sigma$ is a smooth. We denote by $\hat{\cal P}$   the universal
$[S^1]^r$-bundle 
$$\hat{\cal P}:={\cal A}^{**}\times_{{\cal G}}\hat P
$$
over ${\cal B}^{**}\times Y$, and by $\Phi:\hat{\cal P}\ra\C^r$ is the
universal section in the   associated bundle $\hat{\cal
P}\times_{[S^1]^r}\C^r$. For every
$\gamma\in\Z[c_1,\dots,c_r]$ and $h\in H_*(Y)$,  we put as above
$$\delta^\gamma(h)=\Phi^*(\gamma)/h\ .
$$

The image of
$H^*(BK_0,\Z)$ via the natural morphism
$$H^*(BK_0,\Z)\ra H^*(B\hat K,\Z)\ra H^*_{\hat K} (F,\Z)
$$
is the symmetric algebra $S^*(V^*(\Z_m))$.  Using
Proposition 1.1 in [OT2], we see that  that
$\delta^\gamma(h)=0$ if
$\deg(h)<\deg(\gamma)$ and
$\gamma\in S^*(V^*(\Z_m))$.  Therefore the assignment $(\gamma,h)\ra
\delta^\gamma(h)$ descends to a
morphism 
$$\delta:S^*(H)\otimes \Lambda^*(H\otimes H_1(Y,\Z))\ra  H^*({\cal B}^{**},\Z)
$$
of graded algebras, where $H:=\qmod{\   \Z_r}{V^*(\Z_m)}$.

Our invariant is the map
$$GGW(\C^r,\alpha_{\rm can},K_w):[S^*(H)\otimes \Lambda^*(H\otimes
H_1(Y,\Z))]_{2e}\ra\Z\ ,
$$
$$u\mapsto \langle\delta(u),[{\cal M}_{(\tg,g,A^0)}(\lambda)]^{\rm
vir}\rangle\ ,$$
where $e:=\sum_j\chi(L_j)-(r-m)\chi({\cal O}_Y)$ is the expected complex
 dimension of the moduli space, and  $[{\cal
M}_{(\tg,g,A^0)}(\lambda)]^{\rm vir}$ is the   virtual fundamental
class of the zero 
locus of the Fredholm section
$(\mg_\tg,\hg)$. 
\\
\\
\ub{Case} 2.  The action of  $K_w$ on $U(\Sigma)$ is  not free. \\ 

In this case  
let
${\cal G}_{0}\subset {\cal G}$ be the reduced gauge group in a fixed point
$y_0\in Y$, i. e. the kernel of the evaluation map ${\rm ev}_{y_0}:{\cal
G}\ra K_w$, and let
${\cal B}_0:={\cal A}/{\cal G}_0$ be the corresponding quotient.

The same construction  as    above gives a morphism 
$$\delta:S^*(H)\otimes \Lambda^*(H\otimes H_1(Y,\Z))\ra  H^*_{K_w}({\cal
B}_0,\Z)\ ,
$$
and one has a natural restriction map $\rho:H^*_{K_w}({\cal
B}_0,\Q)\ra H^*_{\rm orb}({\cal B}^*,\Q)$. On the other hand, one can generalize
Brussee's method in the orbifold framework and  construct a virtual fundamental
class 
$[{\cal M}_{(\tg,g,A^0)}(\lambda)]^{\rm vir}\in H_{2e}^{\rm orb}({\cal
M}_{(\tg,g,A^0)}(\lambda),\Q)$.

In this case, our invariant is the map
$$GGW(\C^r,\alpha_{\rm can},K_w):[S^*(H)\otimes \Lambda^*(H\otimes
H_1(Y,\Z))]_{2e}\ra\Q\ ,
$$
$$u\mapsto \langle \rho\circ\delta(u),[{\cal M}_{(\tg,g,A^0)}(\lambda)]^{\rm
vir}\rangle\ .$$

\section{A Kobayashi-Hitchin correspondence}

In this section we introduce and study a complex geometric version of the moduli
problem above, and we  prove a Kobayashi-Hitchin type correspondence
which relates the gauge theoretical and the complex geometric moduli spaces.

For $i\in\{1,\dots,m\}$, let ${\cal L}_i^0=(L_i^0,\delta_i^0)$ be a
holomorphic line bundle on $Y$, where $\delta_i^0$ is a fixed semiconnection  on
the differentiable line bundle $L_i^0$.  Let also $(L_j)_{1\leq j\leq r}$ be a  
system of differentiable complex line bundles, and $\lambda$ a system of     
complex isomorphisms
$$\lambda_i:\otimes_{j=1}^r [L_j^{\otimes v^i_j}]\map L_i^0\ .
$$
Our complex geometric moduli problem asks:

Classify all  systems $((\delta_j)_{1\leq i\leq r},(\varphi_j)_{1\leq i\leq r})$, where
$\delta_j$ is a semiconnection in $L_j$ such that
$$\lambda_i(\otimes_{j=1}^r [\delta_j^{\otimes v^i_j}])=\delta_i^0
$$
and $\varphi_j$ is a $\delta_j$-holomorphic section in $L_j$.  Two systems are
considered equivalent if they belong to the same orbit of the complex gauge
group ${\cal G}^\C:={\cal C}^{\infty}(Y,K_w^\C)$, i. e. of

$${\cal G}^\C=\{(f_1,\dots f_r)\in{\cal C}^\infty(Y,\C^*)^r\ |\
\prod\limits_{j=1}^r  f_j^{v^i_j}=1\ \forall i\in\{1,\dots,m\}\}\ .
$$
The configuration space is the product
$$\bar{\cal A}:=\left[\prod_{j=1}^r\bar {\cal
A}(L_j)\right]^{(\lambda,v)}_{\delta^0}\times\bigoplus_{j=1}^r A^0(L_j)
$$
and the equation we consider is just the holomorphy condition
$$\delta_j \varphi_j=0 \ ,\ 1\leq j\leq r\ .\eqno{( H)}
$$

\begin{dt} A system $((\delta_j)_{1\leq i\leq r},(\varphi_j)_{1\leq
i\leq r})$ is called {\it simple} ({\it strictly simple}) if one of the
following equivalent conditions is satisfied:\\ i) its stabilizer  is
finite (respectively trivial),\\ ii)   There is point $y\in Y$ such that
the stabilizer of
$(\varphi_1(y),\dots,\varphi_r(y))$ with respect to the natural
action of $K_w^\C$ in the vector space $\oplus_{j=1}^r L_{j,y}$   is
finite (trivial),\\ 
 iii)   There is point $y\in Y$ such that the stabilizer of
$(\varphi_1(y),\dots,\varphi_r(y))$ with respect to the natural
action of $K_w$ in the vector space $\oplus_{j=1}^r L_{j,y}$  is finite
(trivial).
\end{dt}

\begin{re} A system $((\delta_j)_{1\leq i\leq r},(\varphi_j)_{1\leq
i\leq r})$  is simple if and only if
$$\{\alpha\in\ker(v\otimes\id_\C)|\ 
\alpha_j\varphi_j=0\
\forall j\in\{ 1,\dots, r\}\} =0\ .$$
\end{re}
The simple (strictly simple) systems form an open subspace of $\bar{\cal
A}$ which we denote  by $\bar {\cal A}^{\rm simple}$ ($\bar {\cal A}^{\rm
ssimple}$).

An important role will be played by the moduli spaces  
$${\cal M}^{\rm simple}_{\delta^0}(\lambda):=\qmod{[\bar {\cal A}^{\rm
simple}]^H}{{\cal G}^\C}\ ,\ {\cal M}^{\rm
ssimple}_{\delta^0}(\lambda):=\qmod{[\bar {\cal A}^{\rm
ssimple}]^H}{{\cal G}^\C}
$$
where $[\bar {\cal A}^{\rm
simple}]^H$ stands for the space of simple solutions of the
integrability equation $(H)$.  Note that the data 
$V$, $L_j$, $L_i^0$ can be deduced from $\lambda$, so the notation
${\cal M}^{\rm simple}_{\delta^0}(\lambda)$ makes sense.  

Using standard gauge theoretical methods (see [LO], [OT1], [Su]), one
obtains
\begin{thry} The moduli space ${\cal M}^{\rm
simple}_{\delta^0}(\lambda)$ (${\cal M}^{\rm
ssimple}_{\delta^0}(\lambda)$) has a natural structure of a   complex
analytic orbifold (complex analytic space).
\end{thry}

Note that this complex orbifold (complex space) is in general not smooth
and not Hausdorff.  The local models of the orbifold structure are
(possibly singular) complex spaces endowed with finite group actions; they are
obtained using  {\it local slices} for the ${\cal G}^\C$-action.

There is also an {\it abstract (functorial)} formulation   of our  complex
geometric classification problem, which does not use gauge theoretical methods,
but only classical deformation theory:

Fix an integer matrix $V\in M_{m,r}(\Z)$ as above, let $c=(c_j)_{1\leq j\leq r}$ be
a system of integers, and let ${\cal L}^0=({\cal L}_i^0)_{1\leq i\leq m}$ be a
system of holomorphic line bundles.  

\begin{dt} A holomorphic system of type
$(V,c,{\cal L}^0)$  is  a system 
$$(({\cal L}_j)_{1\leq j\leq r},(\varepsilon_i)_{1\leq i\leq m},(\varphi_j)_{1\leq
j\leq r})\ ,$$
where ${\cal L}_j$ is a holomorphic line bundle of degree $c_j$ on $Y$,
$\varepsilon_i: \bigotimes\limits_{j=1}^r [{\cal L}_j^{\otimes v^i_j}]\ra
{\cal L}_i^0$ is a holomorphic isomorphism, and $\varphi_j\in H^0({\cal
L}_j)$.

An isomorphism between two such systems 
$$(({\cal L}_j)_{1\leq j\leq r},(\varepsilon_i)_{1\leq i\leq
m},(\varphi_j)_{1\leq j\leq r})\ ,\ (({\cal L}'_j)_{1\leq j\leq
r},(\varepsilon_i')_{1\leq i\leq m},(\varphi_j')_{1\leq j\leq r})$$
is  a system of holomorphic isomorphisms
$(u_j)_j$, $u_j:{\cal L}_j\ra{\cal L}'_j$ such that $\varphi'_j=u_j(\varphi_j)$ and
$\varepsilon'_i\circ [\otimes_j  u_j ^{\otimes v^i_j}]=\varepsilon'$.

A system of type $(V,c,{\cal L}^0)$ is called simple (strictly simple) if
its group of automorphisms is finite (trivial).
\end{dt}

Using standards techniques  of deformation theory, one introduces the notion of
families of simple systems of type $(V,c,{\cal L}^0)$ parameterized by a complex
space $Z$ and the notion of isomorphism of such  families.  The functor which
associates to every complex space $Z$ the set of isomorphism classes of families
of simple (strictly simple) systems parameterized by $Z$, is represented
by a complex space ${\cal M}^{\rm simple}_{{\cal L}^0}(V,c)$ (${\cal M}^{\rm
ssimple}_{{\cal L}^0}(V,c)$).

Choosing $c_i=\deg(L_i)$ and ${\cal L}_i^0=(L_i^0,\delta_i^0)$, one gets a
natural embedding  ${\cal M}^{\rm simple}_{\delta^0}(\lambda)\hookrightarrow
{\cal M}^{\rm simple}_{{\cal L}_0}(V,c)$. Note, however, that this embedding is
{\it not surjective}:   the gauge group ${\cal C}^\infty(Y,\C^*)^r$
of $(L_j)_j$ acts in a   natural way on the space $\Lambda$ of systems
$(\lambda_i)_i$ of   complex isomorphisms
$\lambda_i:\bigotimes\limits_{j=1}^r[ L_j^{\otimes v^i_j}]\ra L_i^0$, and one
can easily see that   there is a canonical isomorphism:
$$ {\cal
M}^{\rm simple}_{{\cal L}_0}(V,c)=\coprod_{[\lambda]\in 
\qmod{\scriptstyle\Lambda}{\scriptstyle{\cal C}^\infty(Y,\C^*)^r} }{\cal M}^{\rm
simple}_{\delta^0} (\lambda) \ . 
$$
The quotient $\qmod{\Lambda}{{\cal C}^\infty(Y,\C^*)^r}$ is a
$H^1(Y,\pi_0(K_w))$-torsor, so it has $|\pi_0(K_w)|^{2g(Y)}$ elements.\\

Next we introduce some notations. Let $\Sigma$ be a complete
simplicial fan with $\Sigma(1)\subset\{\R_{\geq 0}v^1,\dots,\R_{\geq
0}v^r\}$, and let $T=(T_j)_{1\leq j\leq r}$ be a system of $r$ complex vector
spaces.

We put
$$U(\Sigma,T):=\{\tau\in\oplus_{j=1}^r T_j|\ 
\exists\sigma\in\Sigma \
\ {\rm such\ that}\ \tau_j\ne0\
\forall   j\in\{1,\dots,r\}\ {\rm for\ which} $$
$$ \R_{\geq 0}v^j \ {\rm
is\ not\ a\ face\ of}\ \sigma  \}\ .
$$
\begin{dt}  Let $\Sigma$ be a complete
simplicial fan in $\R^m$ such that  
$$\Sigma(1)\subset \{\R_{\geq
0}v^1,\dots,\R_{\geq 0}v^r\}\ .$$ 
A system  $(({\cal L}_j)_{1\leq j\leq r},(\varepsilon_i)_{1\leq i\leq
m},(\varphi_j)_{1\leq j\leq r})$ of type $(V,c,{\cal L}^0)$ will be called
$\Sigma$-stable   if one of the following equivalent conditions is
satisfied:
\begin{enumerate}
\item There exists a non-empty Zariski open  set
$Y_0\subset Y$ such that for every $y\in Y_0$ one has 
$(\varphi_1(y),\dots,\varphi_r(y))\in U(\Sigma,{\cal L}_y)$. 
\item $\varphi\in U(\Sigma,H^0({\cal L}))$.
\end{enumerate}
\end{dt}

Here we   denote by $H^0({\cal L})$ the system  $(H^0({\cal L}_j))_{1\leq
j\leq r}$.
The stability condition is obviously an open condition, hence it
defines  open subspaces ${\cal M}^{\Sigma{\raisebox{0.0ex}{-}}\rm
st}_{{\cal L}_0}(V,c)\subset {\cal M}^{\rm simple}_{{\cal L}_0}(V,c)$,  ${\cal
M}^{\Sigma{\raisebox{0.0ex}{-}}\rm st}_{\delta^0}(\lambda)\subset {\cal M}^{\rm simple}
_{\delta_0}(\lambda)$.
\begin{thry} Let $\lambda=\left(\lambda_i:\otimes_{j=1}^r [L_j^{\otimes
v^i_j}]\map L_i^0\right)_i$ be a system of unitary isomorphisms,
where $L_j$, $L_i^0$ are Hermitian line bundles, and
$V=(v^i_j)\in M_{m,r}(\Z)$ is an integer matrix with the properties {\bf
P}$_1$, {\bf P}$_2$. Let $A^0=(A_i^0)_i$,
$A^0_i\in {\cal A}(L_i^0)$ be a system of fixed  Hermitian connections.
Let
$\Sigma$ be a complete simplicial fan with
$\Sigma(1)\subset\{\R_{\geq 0}v^1,\dots,\R_{\geq 0}v^r\}$, and let $\tg\in  K_0(\Sigma)$.
Then there is a natural isomorphism of real analytic orbifolds
$${\cal M}_{(\tg,g,A^0)}(\lambda)\simeq {\cal M}^{\Sigma{\raisebox{0.0ex}{-}}\rm
st}_{\bar\partial_{A^0}}(\lambda)\ .
$$
 \end{thry}
\pf  By the universal Kobayashi-Hitchin correspondence for universal
vortices [Mu1], [LT2], there is a natural isomorphism of real analytic
orbifolds
$${\cal M}_{(\tg,g,A^0)}(\lambda)\simeq {\cal M}^{\rm
st}_{\mu_w+i\tg,\bar\partial_{A^0}}(\lambda)\ ,
$$
where the right hand side denotes  the moduli space of integrable pairs
$(\delta,\varphi)$ which are analytically stable with respect to the moment
map
$\mu_w+i\tg$, $z\mapsto -\frac{i}{2}p_v( |z_1|^2,\dots,|z_r|^2)+i\tg$.

Recall that a pair $(\delta,\varphi)$ is $(\mu_w+i\tg)$-analytically stable  
if for every $\xi\in \ker v$ one has 
$$\int\limits_Y
\lambda_{\mu_w+i\tg}(\varphi,
s(\xi))vol_g>0\ .
$$
But 
$$\lambda_{\mu_w+i\tg}(\varphi,
s(\xi))(y)=\lim_{t\ra\infty}\left\langle
p_v\left(\matrix{\frac{1}{2}e^{2\xi_1 t}|\varphi_1(y)|^2\cr
\cdot\cr\cdot\cr\cdot\cr\frac{1}{2}e^{2\xi_r t}
|\varphi_r(y)|^2}\right)-\tg,\xi\right\rangle=$$
$$=\left\{\begin{array}{cl}
-\langle\tg,\xi\rangle &{\rm if}\ \xi_j\leq 0\ \forall j \ 1\leq j\leq
r\ {\rm for\ which}\ \varphi_j(y)\ne 0\\ \infty  &{\rm if}\ \exists   j
\ 1\leq j\leq r\ {\rm such\ that}\ \varphi_j(y)\ne 0\ {\rm and}\ \xi_j>0\ .
\end{array}\right.
$$

This shows that $(\delta,\varphi)$ is $(\mu_w+i\tg)$-analytically
stable if and only if for every $\xi\in \ker v$ for which
$\langle\tg,\xi\rangle\geq 0$ there exists a point $y\in Y$ and an index
$j\in\{1,\dots,r\}$ such  that $\varphi_j(y)\ne 0$   and $ \xi_j>0$.  It is
easy to see, using the holomorphy of $\varphi$, that this happens if and only
if there exists a point
$y\in Y$ such that   for every $\xi\in \ker v$ for which
$\langle\tg,\xi\rangle\geq 0$ there is some  $j\in\{1,\dots,r\}$ such  that
$\varphi_j(y)\ne 0$   and $ \xi_j>0$. This holds if and only if there exists
$y\in Y$ such that $\varphi(y)$ is analytically stable with respect to the
natural action of $K_w$ on
$\oplus_{j=1}^r L_{j,y}$ and the same moment map $\mu_w+i\tg$.
\qed

\section{Complex geometric description of the moduli spaces}

\subsection{Linear spaces associated with Poincar\'e line bundles}

Let ${\cal E}$ be a holomorphic bundle over the product 
$X\times P$ of two complex manifolds, with $X$ compact and finite dimensional.
The disjoint union
$$H^0_X({\cal E}):=\coprod_{p\in P} H^0({\cal E}|_{\{p\}\times X})
$$
has a natural structure of a holomorphic  linear space over $P$.  There
are two ways to construct this linear space:
\begin{enumerate}
\item Using infinite dimensional analytic geometry:

We define $H^0_X({\cal E})$ as the vanishing locus of the relative
$\bar\partial$-equation equation as follows:

Fix  Hermitian structures on $X$ and ${\cal E}$. For every $(p,q)$ and
$k>0$ the union
$A^{pq}_X({\cal E})_k:=\coprod_{p\in P} L^2_k({\cal E}|_{\{p\}\times
X}\otimes\Lambda^{pq}_X)$ is naturally a holomorphic locally trivial Banach
bundle over $P$.  The holomorphic sections in $A^{pq}_X({\cal E})_k$ over an
open set
$Q\subset P$ are the sections in $[{\rm pr}_X^*(\Lambda^{pq}_X)\otimes{\cal
E}]|_{Q\times X}$ which are $L^2_k$-Sobolev in the $X$-direction
and holomorphic in the $Q$ direction.

Then $H^0_X({\cal E})$ is defined as the vanishing locus  of the relative
differential 
$$\bar\partial_{{\cal E}/X}:=A^{00}_X({\cal E})_k\ra
A^{01}_X({\cal E})_{k-1}$$
  The result does not depend on
the Sobolev index $k\in \N$.

\item  Using the duality between linear spaces and coherent sheaves:
 
One defines  $H^0_X({\cal E})$ as the linear space which corresponds to
the sheaf 
$$R^{\dim(X)} [{\rm pr}_P]_*({\cal E}^{\vee}\otimes{\rm pr}_X^*(K_X))\
.$$

Note that  if  $H^0({\cal E}|_{\{p\}\times X})=0$ for generic $p\in P$,
then   $[{\rm pr}_P]_*({\cal E})=0$, whereas the linear space
$H^0_X({\cal E})$ is non-trivial as soon $H^0({\cal E}|_{\{p\}\times
X})\ne 0$ for at least one
$p\in P$.

The equivalence between the two definitions follows from the explicit
description of the structure sheaf of the complex subspace defined by
an equation with values in an infinite dimensional complex Banach space
(see [Dou],   [LT]). 

\end{enumerate}

Let now $Y$ be a complex curve, and fix a point $y_0\in Y$. For every
class $c\in H^2(Y,\Z)$ denote by
$\Pg^c_{y_0}$ the Poincar\'e line bundle over $Y\times\Pic^c(Y)$ which is
trivial on $\{y_0\}\times \Pic^c(Y)$. Denote by $\Lg^c_{y_0}$ the
corresponding linear space over $\Pic^c(Y)$, i. e.
$$\Lg^c_{y_0}=H^0_Y(\Pg^c_{y_0}).
$$

The Poincar\'e line bundle $\Pg^c_{y_0}$ has a gauge theoretical
description, which we sketch here briefly:

Let $L^c$ be a differentiable line
bundle of Chern class $c$, and let $\bar{\cal A}(L^c)_k$ be the $L^2_k$
completion of the affine space of semiconnections on $L^c$  for a
sufficiently large Sobolev index $k$. Denote by ${\cal
C}_{y_0}(Y,\C^*)_{k+1}\subset {\cal C}(Y,\C^*)_{k+1}$ 
the kernel of the evaluation morphism
${\rm ev}_{y_0}$ on the $(k+1)$-Sobolev completion ${\cal
C}(Y,\C^*)_{k+1}$ of the gauge group ${\cal C}^\infty(Y,\C^*)$.  The
group ${\cal
C}_{y_0}(Y,\C^*)_{k+1}$   acts freely on  
$\bar{\cal A}(L^c)_k$.

The line bundle
${\cal T}:={\rm pr}_Y^*(L^c)$ on the product  $Y\times\bar{\cal
A}(L^c)_k$ comes with  a tautological integrable semiconnection, which
agrees with
$\delta$ over $Y\times\{\delta\}$ for every   $\delta\in\bar{\cal A}(L_c)$,
and which is
  the standard  
$\bar\partial$-operator over the infinite dimensional fibres $\{y\}\times
\bar{\cal A}(L^c)_k$.  

The free action of ${\cal
C}_{y_0}(Y,\C^*)_{k+1}$ on $Y\times\bar{\cal A}(L^c)_k$ admits a
tautological holomorphic linearization in the line bundle ${\cal T}$.
  Therefore,  as in the finite dimensional framework, one can
regard the ${\cal
C}_{y_0}(Y,\C^*)_{k+1}$-quotient of ${\cal T}$ as a holomorphic
line bundle over 
$$\qmod{Y\times\bar{\cal
A}(L^c)_k}{{\cal C}_{y_0}(Y,\C^*)_{k+1}}=Y\times\Pic^c(Y)\ .$$
This quotient line bundle  is just the Poincar\'e bundle
$\Pg^c_{y_0}$.
\begin{re} Using this gauge theoretical description of the Poincar\'e
line bundle, one gets  a corresponding gauge  theoretical description
of the associated linear space $\Lg^c_{y_0}$. This linear space can be
identified with the finite dimensional subspace cut out by the
holomorphy condition
$\delta\varphi=0$ in the quotient
$$\qmod{\bar{\cal A}(L^c)_k\times A^0(L^c)_k}{{\cal
C}_{y_0}(Y,\C^*)_{k+1}}\ .
$$
\end{re}

Consider now the divisor $P^c_y:=\{y\}\times\Pic^c(Y)$ of $Y\times\Pic^c(Y)$.
Any point
$y\in Y$ defines a topologically trivial holomorphic line bundle 
$$\Pg^{c}_{y_0y}:=
\Pg^{c}_{y_0}|_{ P^c_y}$$
over $\Pic^c(Y)$. Let $\pi^c:Y\times\Pic^c(Y)\ra\Pic^c(Y)$ be the
natural projection.

We have the exact sequence of sheaves on $Y\times\Pic^c(Y)$:
$$0\map[\Pg^{c}_{y_0y}]^\vee\otimes{\rm pr}_Y^*(K_Y)\map
[\Pg^{c}_{y_0y}]^\vee\otimes{\rm pr}_Y^*(K_Y) (P^c_y)\map$$
$$\map[\Pg^{c}_{y_0y}]^\vee\otimes{\rm pr}_Y^*(K_Y)(P^c_y)_{P^c_y}\map 0\ ,
$$
and  a canonical isomorphism 
$${\rm
pr}_Y^*(K_Y)(P^c_y)_{P^c_y}={\cal O}_{P^c_y}\ ,$$
 which is induced by the obvious
isomorphism
$[K_Y](\{y\})_{\{y\}} \simeq\C$.  Therefore we get a
morphism
$$[\Pg^{c}_{y_0y}]^\vee=
\pi^c_*([{\Pg^c_{y_0}}]^\vee_{P^c_y}) 
\map R^1\pi^c_*([{\Pg^c_{y_0}}]^\vee\otimes ({\rm pr}_Y)^*K_Y)\ ,$$
which induces a morphism of linear spaces  
$${\rm
ev}^{c}_{y_0y}:\Lg^c_{y_0}\ra \Pg^{c}_{y_0y}$$
over $\Pic^c(Y)$. As a map of sets, ${\rm
ev}^{c}_{y_0y}$ is just the evaluation map  associated with the point $y$.
\begin{pr} Let $\iota_y:\Pic^{c}(Y)\ra\Pic^{c+1}(Y)$ be the isomorphism
defined by $\iota_y([{\cal L}]):=[{\cal L}(\{y\})]$. 
\begin{enumerate}
\item There is a canonical isomorphism 
$$(\id_Y\times\iota_y)^*(\Pg_{y_0}^{c+1}(-P^{c+1}_y)) \simeq \Pg_{y_0}^{c}\ .
$$
\item There is an exact sequence of linear spaces over $\Pic^c(Y)$
$$0\map\Lg^c_{y_0}\map\iota_y^*(\Lg^{c+1}_{y_0})\textmap{\iota_y^*({\rm
ev}^{c+1}_{y_0y})}\iota_y^*(\Pg^{c+1}_{y_0y})\ . 
$$
\end{enumerate}
\end{pr}
\pf The first property follows from the universal property of the Poincar\'e line
bundle.  The second property can be obtained as follows:

Using the   exact sequence of line bundles over $Y\times\Pic^{c+1}(Y)$
$$0\map  [\Pg_{y_0}^{c+1}]^\vee\otimes {\rm
pr}_Y^*(K_Y)\map [\Pg_{y_0}^{c+1}]^\vee\otimes {\rm
pr}_Y^*(K_Y)(P^{c+1}_y)\map  $$
$$\map [\Pg_{y_0}^{c+1}]^\vee\otimes
{\rm pr}_Y^*(K_Y)(P^{c+1}_y)_{P^{c+1}_y} \map 0
$$
one gets an exact sequence of sheaves on $\Pic^{c+1}(Y)$:
$$ \pi^{c+1}_*([\Pg_{y_0}^{c+1}]^\vee\otimes {\rm
pr}_Y^*(K_Y)(P^{c+1}_y))\map
 \pi^{c+1}_*([\Pg_{y_0}^{c+1}]^\vee\otimes {\rm
pr}_Y^*(K_Y)(P^{c+1}_y)_{P^{c+1}_y})\map
$$
$$\map R^1\pi^{c+1}_*([\Pg_{y_0}^{c+1}]^\vee \otimes {\rm
pr}_Y^*(K_Y))\map R^1\pi^{c+1}_*([\Pg_{y_0}^{c+1}]^\vee\otimes {\rm
pr}_Y^*(K_Y)(P^{c+1}_y))\map 0\ .
$$
The claim follows now from the first statement and   the
definition of the morphism
${\rm ev}^{c+1}_{y_0y}$.
\qed
\begin{co} Let $\{y_1,\dots,y_k\}$ be a finite set of $Y$ and let $\iota_D$ be the
 isomorphism $\Pic^c(Y)\ra\Pic^{c+k}(Y)$ defined by the divisor $D=\sum_i y_i$.
Then there is a canonical isomorphism  
$$(\iota_D)_*(\Lg^c_{y_0})\simeq\ker\left[\bigoplus\limits_{i=1}^k{\rm
ev}^{c+k}_{y_0y_i}:\Lg^{c+k}_{y_0}\map\bigoplus\limits_{i=1}^k \Pg^{c+k}_{y_0y_i}
\right]\ .
$$
\end{co} 
Note that the linear space $\Lg^{c+k}_{y_0}$ becomes a vector bundle if
$c+k>2(g-1)$. This remark will play an important role in the following
section.

\subsection{Moduli spaces of stable systems as subspaces of locally
trivial toric fibre bundles}

Let $B$ be  a complex space, let ${\cal T}=({\cal T}_j)_{1\leq
j\leq r}$ be a system of holomorphic linear spaces over $B$, and let $\Sigma$
be a complete simplicial fan with  
$$\Sigma(1)\subset\{\R_{\geq
0}v_1,\dots,\R_{\geq 0}v_r\}\ .$$

Note that the union $U_B(\Sigma,{\cal T}):=\union_{b\in B} U(\Sigma,{\cal
T}_b)$ is an open subspace of the fibre product
${\prod\limits_{j=1}^r}{}_B{\cal T}_j$. Since the natural $K_w^\C$-action on
$U_B(\Sigma,{\cal T})$ admits local slices,  the quotient 
$$X_{\Sigma,B}({\cal T}):=\qmod{U_B(\Sigma,{\cal T})}{K^\C_w}
$$
has a  natural 
complex analytic (in
general singular) orbifold structure. It comes with a natural proper morphism to
$B$.

\begin{dt}   The morphism $X_{\Sigma,B}({\cal T}) \map B$
will be called the toric
fibration associated with the data $(w,\Sigma,{\cal T})$.
\end{dt}

Consider the map
$$\Pic(V):\prod_{j=1}^r\Pic^{c_j}(Y)\map\prod_{i=1}^m\Pic^{c_1(L_i^0)}(Y)
$$
defined by the matrix $V$.  The complex manifold
$$P:=\qmod{\left[\prod_{j=1}^r\bar {\cal
A}(L_j)\right]^{(\lambda,v)}_{\delta^0}}{{\cal G}^\C}
$$
can be identified with a connected component  of the fibre 
$\Pic(V)^{-1}([L^0,\delta^0])$ of
$\Pic(V)$; it is an abelian variety of dimension $g(Y)(r-m)$.

Denote by $\underline c$ the
system $(c_j)_{1\leq j\leq r}$,  by ${\Lg}^{\underline c}$ the system
$({\Lg}^{c_j}_{y_0}))_{1\leq j\leq r}$, and by $p_{\underline c}$ the system
of projections $p_{c_j}:\prod_{j=1}^r\Pic^{c_j}(Y)\ra \Pic^{c_j}(Y)$.

\begin{thry}
Choose $\tg\in K_0(\Sigma)$ and fix $y_0\in Y$.  There is a  canonical
isomorphisms of complex analytic orbifolds 

$${\cal M}^{\Sigma\hbox{-}\rm
st}_{\delta^0}(\lambda)\simeq
X_{\Sigma,P}( p_{\underline c}^*({\Lg}^{{\underline c}}_{y_0})|_P)\ .$$
\end{thry}
\pf We denote by ${\cal G}^\C_{y_0}\subset{\cal G}^\C$ the kernel of the
evaluation morphism $ev_{y_0}:{\cal G}^\C\ra K_w^\C$ associated with   $y_0$.
${\cal G}^\C_{y_0}$ acts freely on the configuration space 
$\left[\prod_{j=1}^r\bar {\cal
A}(L_j)\right]^{(\lambda,v)}_{\delta^0}\times\oplus_{j=1}^r A^0(L_j)$.

  Denote by $\left[\prod_{j=1}^r\bar {\cal
A}(L_j)\right]_{P}$ the space of   systems of semiconnections
$\delta=(\delta_j)_j\in \prod_{j=1}^r\bar {\cal
A}(L_j) $ such that $([\delta_j])_j\in P$. The larger configuration
space
$$\left[\prod_{j=1}^r\bar {\cal
A}(L_j)\right]_{P}\times\bigoplus_{j=1}^r A^0(L_j)
$$
comes with a free action of ${\cal C}^\infty_{y_0}(Y,\C^*)^{ r}$, and
using the gauge theoretical interpretation of the linear spaces associated
with Poincar\'e line bundles (Remark 4.1) , one sees easily that the
complex subspace cut out by the integrability condition
$\delta\varphi=0$ in the quotient
$$\bigqmod{\left[\prod_{j=1}^r\bar {\cal
A}(L_j)\right]_{P}\times\bigoplus_{j=1}^r A^0(L_j)}{{\cal
C}^\infty_{y_0}(Y,\C^*)^{r}}
$$
is precisely
$$ \left.\left[\prod_{j=1}^r {\Lg}^{c_j}_{y_0}\right] \right|_P
=\prod\limits_{j=1}^r{}_{_{\scriptstyle P}}
((p_{c_j}^*({\Lg}^{c_j}_{y_0})|_P)_j)\ .
$$
On the other hand, after suitable Sobolev completions,  the natural morphism
$$\bigqmod{\left[\prod_{j=1}^r\bar {\cal
A}(L_j)\right]^{(\lambda,v)}_{\delta^0}\times\bigoplus_{j=1}^r
A^0(L_j)}{{\cal G}^\C_{y_0}}\ra\bigqmod{\left[\prod_{j=1}^r\bar {\cal
A}(L_j)\right]_{P}\times\bigoplus_{j=1}^r A^0(L_j)}{{\cal
C}^\infty_{y_0}(Y,\C^*)^{r}}
$$
becomes an isomorphism of Banach manifolds, and induces an isomorphism of
complex spaces
$$\left\{[\delta,\varphi]\in\left.\bigqmod{\left[\prod_{j=1}^r\bar {\cal
A}(L_j)\right]^{(\lambda,v)}_{\delta^0}\times\bigoplus_{j=1}^r
A^0(L_j)}{{\cal G}^\C_{y_0}}\right|\
\delta\varphi=0\right\}\textmap{\simeq}
\prod\limits_{j=1}^r{}_{_{\scriptstyle P}}
((p_{c_j}^*({\Lg}^{c_j}_{y_0})|_{P})_j)\ .
$$

According to Theorem 3.6,  a class $[\delta,\varphi]$ is
$\mu_w+i\tg$-stable if and only if the corresponding element in this
fibre product belongs to $U_P(\Sigma,p_{\underline
c}^*({\Lg}^{\underline c}_{y_0})|_P )$.  Therefore one obtains an
isomorphism
$${\cal M}^{\Sigma\hbox{-}\rm
st}_{\delta^0}(\lambda)\simeq
\qmod{U_P(\Sigma,p_{\underline c}^*({\Lg}^{\underline c}_{y_0})|_P)}{K^\C_w}\ .
$$
\qed

Now fix  effective divisors $D_j$ on $Y$ of sufficiently large degrees $d_j$  
such that $h^1({\cal L'})=0$ for every holomorphic line bundle ${\cal L}'$ of
degree $c'_j:=c_j+d_j$. We   assume that $D_j$ is   a set of $d_j$
distinct simple points $y_j^i$, $1\leq i\leq d_j$.  Denote by $\underline D$ the
system of divisors
$D_j$ and by 
$$\iota_{\underline D}:\prod_{j=1}^r\Pic^{c_j}(Y)\ra
\prod_{j=1}^r\Pic^{c'_j}(Y)$$
the isomorphism  defined by the system  of maps  $\otimes{\cal
O}(\underline D)$.  Set $P':=\iota_{\underline D}(P)$.

The  toric fibration $p':X_{\Sigma,P'}( p_{\underline c'}^*({\Lg}^{{\underline
c'}}_{y_0})|_{P'})\ra P'$ is a locally trivial fibre bundle over $P'$. Let
$q_j$ be the natural map  
$$q_j:{\prod_s} \Lg^{c'_s}_{y_0}\ra \Pic^{c'_j}(Y)\ .$$

 The maps 
${\rm ev}^{c_j'}_{y_0,y_j^i}:{\Lg}^{c'_j}_{y_0}\ra \Pg^{c'_j}_{y_0 y^i_j}$ are
equivariant with respect to the natural actions of the group $K^\C_w$, so they
define sections $e^i_j$ in the orbifold line bundles
$${\cal L}^i_j:=\bigqmod{\left.[q_j]^*(\Pg^{c'_j}_{y_0
y^i_j})\right|_{U_{P'}(\Sigma,p_{\underline c'}^*({\Lg}^{\underline
c'}_{y_0}) ) }} {K^\C_w} 
$$
over the {\it locally trivial} toric fibre bundle $X_{\Sigma,P'}( p_{\underline c'}^*({\Lg}^{{\underline
c'}}_{y_0})|_{P'})$.
\\

Combining Corollary 4.3 and Theorem 4.5, we get the following description of our
moduli space as the vanishing locus of a system of sections in line bundles over
a locally trivial toric fibre bundle.
\begin{thry} (Embedding theorem) There are canonical isomorphisms
$${\cal M}^{\Sigma\hbox{-}\rm
st}_{\delta^0}(\lambda)\simeq
X_{\Sigma,P}( p_{\underline c}^*({\Lg}^{{\underline c}}_{y_0})|_P)\simeq Z((e^i_j)
_{\matrix{\scriptstyle
1\leq i\leq d_j\vspace{-1mm}\cr\scriptstyle
1\leq j\leq r}} )\ .
$$
\end{thry} 

Using these isomorphisms and the methods developed by Fulton ([Fu], p. 244)
and Behrend-Fantechi, one can endow 
${\cal M}^{\Sigma\hbox{-}\rm st}_{\delta^0}(\lambda)$ with a
distinguished homology class, namely the virtual fundamental class associated
with  the section  
$\oplus_{i,j}(e^i_j)$ in the bundle
$\oplus_{i,j} {\cal L}^i_j$ over the smooth orbifold $X_{\Sigma,P'}(
p_{\underline c'}^*({\Lg}^{{\underline c'}}_{y_0})|_{P'})$.  We will
call this class {\it the algebraic geometric  virtual fundamental
class}. Note that  the construction of the embedding in Theorem 4.6
depends on the (non-canonical) parameters $y_0$ 
 and $D_j$, hence the obtained virtual fundamental class might  a priori
depend also on these parameters.

\section{Identifying virtual fundamental classes}

The purpose of the this chapter is to show that, at least in the smooth
case,  the algebraic geometric virtual fundamental class obtained in the
previous section coincides with the correct virtual fundamental class of 
our moduli space, namely with the virtual fundamental class associated
with  its initial gauge theoretical construction (see section 2.1). 
Consequently, the algebraic geometric virtual fundamental class does not depend
on the non-canonical parameters  $y_0$  
and $D_j$ which occur in the construction of the embedding of the
moduli space in the smooth toric fibre bundle $X_{\Sigma,P'}(
p_{\underline c'}^*({\Lg}^{{\underline c'}}_{y_0})|_{P'})$.\\

We come back to the gauge theoretical construction of the moduli spaces 
${\cal M}^{\Sigma\hbox{-}\rm st}_{\delta^0}(\lambda)$:

The configuration space of our moduli problem is
$$\bar{\cal A}=\left[\prod_{j=1}^r\bar {\cal
A}(L_j)\right]^{(\lambda,v)}_{\delta^0}\times\bigoplus_{j=1}^r
A^0(L_j)=[\bar{\cal A}(L)]^{(\lambda,v)}_{\delta^0}\times A^0(L)\ ,
$$
where  $L:=\oplus_j L_j$ is regarded as a rank $r$ vector bundle with
structure group $[\C^*]^r$.

The space 
$$\left[ A^{pq}(Y)^{\oplus r}\right]^v:=\{(\alpha_1,\dots,\alpha_r)\in 
A^{pq}(Y)^{\oplus r}|\
\sum_j v^i_j \alpha_j=0\}
$$
is the model vector space of the affine space  $[\bar{\cal
A}(L)]^{(\lambda,v)}_{\delta^0}$.

In order to complete our configuration space, we fix
a large integer $k\gg 0$ and a real number $s\in (1,2)$. 

We denote by
$\left[ A^{pq}(Y)^{\oplus r}\right]^v_k$ ($[\bar{\cal
A}(L)]^{(\lambda,v)}_{\delta^0,k}$ ) the completion of the  vector
space
$\left[ A^{pq}(Y)^{\oplus r}\right]^v$ (affine space $[\bar{\cal
A}(L)]^{(\lambda,v)}_{\delta^0}$ ) with respect to the Sobolev norm $L^2_k$
(Sobolev
  $L^2_k$-topology), and by $L^s(L)$ the completion of $A^0(L)$
with respect to the Sobolev norm $L^s$.   

Thus our completed configuration space is
$$\bar{\cal A}^s_k:=[\bar{\cal
A}(L)]^{(\lambda,v)}_{\delta^0,k}\times L^s(L)\ .
$$

The reason for
completing the configuration space in this odd way will become clear
later: the completed configuration space contains in particular pairs whose
section  component is   meromorphic with simple poles.  

Note that, if we
fix a Hermitian metric on $L$, then   the function $|\varphi|^2$
associated to an $L^s$-section $\varphi$ is not integrable, so the vortex
equation is not well defined on the completed configuration space. In
particular,  one {\it cannot} use such a Sobolev completion to prove a
Kobayashi-Hitchin correspondence.

As usual we let ${\cal G}^\C_{k+1}$ be the completion of the  gauge
group with respect to the $L^2_{k+1}$ norm.  

The left hand term $\delta\varphi$ in the
integrability equation 
$$\delta\varphi=0\eqno{(  H)}
$$
should be regarded as an element of  the
distribution space
$L^s_{-1}(\Lambda^{01}_Y\otimes L)$.

One can show that the completion
procedure above does not introduce new orbits in the moduli space, hence
\begin{pr} Any
weak solution $(\delta_1,\varphi_1)\in\bar{\cal A}^s_k$ of the integrability
equation $(H)$  is ${\cal G}^\C_{k+1}$-equivalent to a
smooth  solution, which is unique up to  ${\cal
G}^\C$-equivalence. 
\end{pr}  
\pf To prove this, one  brings $\delta_1$ in  
Coulomb gauge with respect to a smooth semiconnection $\delta_0$ using a
gauge transformation of the type $\exp(u)$, $u\in \left[ A^{00}(Y)^{\oplus
r}\right]^v_{k+1}$.  The new semiconnection $\delta_1':=\delta\cdot\exp(u)$
will be smooth, since the $(0,1)$ form $\delta_1'-\delta_1$ solves an
elliptic equation with smooth coefficients. The $L^s$-section
$\varphi_1':=\varphi_1\cdot\exp(u)$ will  also be smooth, by elliptic
regularity.
 \qed
\begin{co} If $(\delta,\varphi)\in\bar{\cal A}^s_k$ is an \ub{integrable}
pair, then $\varphi\in L^2_{k+1}(L)$.
\end{co}

We suppose for simplicity that the $K^\C_w$-action  on
$U(\Sigma)$ is free, so that all (weak) solutions  of the system
$(H, V_g^\tg)$ have trivial stabilizers.  Therefore they belong to
the subspace $[\bar{\cal A}^s_k]^{**}$ of weak pairs with
trivial stabilizer.

Since we use a very weak Sobolev topology on the second factor of
our configuration space,  the general  
procedure ([DK], [LT]) to endow the quotient
$$\bar{\cal B}^{**}:=\qmod{[\bar{\cal A}^s_k]^{**}}{{\cal
G}^\C_{k+1}}
$$
with the structure of a Banach manifold must be adapted carefully. The
usual (formal)
$L^2$-adjoint of the infinitesimal action 
$$D^0_p(f)=(\bar\partial f,- f\varphi)
$$
at a point $p=(\delta,\varphi)$ is
given by the formula
$$[D^0_p]^*(\alpha,\psi)=\bar\partial^*\alpha-\sum_j (\psi_j,\varphi_j)\ ,
$$
where the pairing $(\cdot,\cdot)$ stands for the pointwise Hermitian
product.  The difficulty comes from the fact that, when
$\varphi$ is only an 
$L^s$-section, $[D^0_p]^*$ does not necessarily extend   to a bounded
operator
$\left[ A^{01}(Y)^{\oplus r}\right]^v_k\times L^s(L)\ra \left[
A^{00}(Y)^{\oplus r}\right]^v_{k-1}$.

However, one can  modify this adjoint operator  to get bounded
operators with the desired properties. The restriction of $D^0_p$ to the Lie
algebra   
$${\rm Lie}(K_w^\C)=\ker v\otimes\C$$
 of the subgroup of constant
gauge transformations is given by 
$$f\mapsto (0,-f\varphi)\ .
$$
The   operator $m_\varphi:\ker v\otimes\C\ra L^s(L)$ defined by  
multiplication with $\varphi$ is injective, since
$(\delta,\varphi)$ is an irreducible pair. The image of this operator is a
finite dimensional subspace of the Banach space $L^s(L)$, hence it is 
closed and has a closed complement. Therefore, $m_\varphi$ admits a
continuous left inverse, say $q_\varphi:L^s(L)\ra \ker v\otimes\C$. We set
$$q_p(\alpha,\psi):= \bar\partial^*\alpha-q_\varphi(\psi)\ , \ q_p:\left[
A^{01}(Y)^{\oplus r}\right]^v_k\oplus L^s(L)\ra \left[ A^{00}(Y)^{\oplus
r}\right]^v_{k-1}\ .
$$
Note that, since the operators $\bar\partial^*$ and $q_\varphi$ take values in
direct summands, one has 
$$\ker q_p=\ker\bar\partial^*\oplus\ker q_\varphi\ .\eqno{(3)}
$$
\begin{pr}  Put
$$V_{p,\varepsilon}:=\{p+(\alpha,\psi)|\ q_p(\alpha,\psi)=0,\
\nr\alpha\nr_{L^2_k}<\varepsilon,\ \nr\psi\nr_{L^s}<\varepsilon\}\ ,
$$
and define $f_{p,\varepsilon}:V_{p,\varepsilon}\times{\cal G}^\C_{k+1}\ra
[\bar{\cal A}^s_k]$ by   the formula
$(p',\gamma)\mapsto p'\cdot
\gamma$.
\begin{enumerate}
\item Let $p=(\delta,\varphi)\in[\bar{\cal A}^s_k]^{**}$ be a weak 
strictly irreducible  pair.  Then there exists $\varepsilon_p>0$ such
that $f_{p,\varepsilon_p}$ defines a diffeomorphism from 
$V_{p,\varepsilon}\times{\cal G}^\C_{k+1}$ onto an open neighbourhood of 
the orbit $p\cdot{\cal G}_{k+1}^\C$ in
$[\bar{\cal A}^s_k]^{**}$.

\item There is a unique Banach manifold structure on $\bar{\cal
B}^{**}$ such that the natural maps  $V_{p,\varepsilon_p}\map\bar{\cal
B}^{**}$ become smooth parametrizations.
\end{enumerate}
\end{pr}
\pf   We  first seek $\varepsilon$ sufficiently small such that
$f_{p,\varepsilon}$ becomes injective.  Let  $(\alpha_i,\psi_i)\in\ker
q_p$,
$g_i\in{\cal G}_{k+1}^\C$ such that  
$$(p+(\alpha_1,\psi_1))\cdot g_1=(p+(\alpha_2,\psi_2))\cdot g_2\ .$$ 
Put $\gamma:=g_1g_2^{-1}$.  It follows that
$$\gamma^{-1}\bar\partial\gamma=\alpha_2-\alpha_1\ ,\ 
(\gamma^{-1}-1)\varphi=\psi_2-\psi_1\ .\eqno{(4)}$$

Write $\gamma=k+\gamma_0$, where $k\in[\C]^r$ and
$\gamma_0\in\bar\partial^*(A^{01}(Y,\C^r)_{k+2})$. The first relation  can
be written as
$$\bar\partial\gamma_0=k
(\alpha_2-\alpha_1)+\gamma_0(\alpha_2-\alpha_1)\ .  \eqno{(5)}
$$
Taking into account that $(\alpha_i,\varphi_i)\in\ker q_p$ and using (3), we
get  $\bar\partial^*(\alpha_2-\alpha_1)=0$, so (5) yields
$$\Delta\gamma_0=\bar\partial^*[\gamma_0(\alpha_2-\alpha_1)]\ .
$$
This gives an estimate of the form
$$\nr \gamma_0\nr_{L^2_{k+1}}\leq
c\nr\gamma_0(\alpha_2-\alpha_1)\nr_{L^2_k}\leq
c'\nr\gamma_0\nr_{L^2_{k+1}}\nr(\alpha_2-\alpha_1)\nr_{L^2_k}\ .
$$
If we choose $\varepsilon<\frac{1}{2c'}$, this inequality implies
$\gamma_0=0$, hence $k=\gamma\in K_w^\C$.  The second relation in (4)
can be written as $m_\varphi(\gamma-1)=\psi_2-\psi_1$.  Using  
 $(\alpha_i,\psi_i)\in\ker q_p$ and (3), one obtains
$$(\gamma-1)=q_\varphi m_\varphi(\gamma-1)=q_\varphi(\psi_2-\psi_1)=0\ .
$$
This shows that $g_1=g_2$, and $(\alpha_1,\psi_1)=(\alpha_2,\psi_2)$.\\

Now we want to prove that, for sufficiently small $\varepsilon$,
$f_{p,\varepsilon}$ is \'etale. It suffices to prove this property  at
the points of the form $(p',e)$. Indeed,  since
$f_{p,\varepsilon}(p',h\gamma)= f_{p,\varepsilon}(p',h)\cdot \gamma$,
it follows that  $f_{p,\varepsilon}$ is \'etale in $(p',\gamma)$ if and
only if  it is \'etale in $(p',e)$.  Write $p'=p+(\alpha',\psi')$. One
has
$$d_{(p',e)}f_{p,\varepsilon}(\alpha,\psi,u)=(\alpha+\bar\partial
u,\psi-u(\varphi+\psi'))\ .
$$
For $\psi'=0$, this operator is invertible. Its inverse 
$$r:\left[A^{01}(Y)^{\oplus r}\right]^v_k\oplus L^s(L)\map\ker q_p\oplus
\left[ A^{00}(Y)^{\oplus r}\right]^v_{k+1}
$$
is given by
$$r(a,v)=\left(a-\bar\partial G\bar\partial a,\left[G\bar\partial a-
q_\varphi\left((G\bar\partial a) \varphi+v\right)\right]\varphi+v
,G\bar\partial a- q_\varphi\left((G\bar\partial a)
\varphi+v\right)\right)\ ,
$$
where $G$ denotes the Green operator of the Laplacian
$\Delta=\bar\partial^*\bar\partial$.

Since the linear bounded operator $d_{(p',e)}f_{p,\varepsilon}$ depends
continuously on the parameter $\psi'\in L^s(L)$, it
must remain invertible for $\psi'$ sufficiently small.

The second statement follows easily from the first.
\qed

The map $\bar h:[\bar{\cal A}^s_k]^{**}\ra L^s_{-1}(\Lambda^{01}_Y\otimes
L)$,  given by $(\delta,\varphi)\mapsto \delta\varphi$,
is equivariant with respect to the action of the gauge group ${\cal
G}^\C_{k+1}$, hence it descends to a holomorphic Fredholm  section
$\bar\hg$ in the Banach bundle
$$\bar \Eg:=[\bar{\cal A}^s_k]^{**}\times_{{\cal
G}^\C_{k+1}}L^s_{-1}(\Lambda^{01}_Y\otimes L)
$$
over $\bar{\cal B}^{**}$.

Using the  regularity result given by Proposition 5.1, it follows that
the  moduli space ${\cal M}^{\rm ssimple}_{\delta^0}(\lambda)$ is  
the vanishing locus of the   holomorphic Fredholm  section $\bar \hg$.
Under our assumptions, the moduli space ${\cal M}^{\Sigma\hbox{-}\rm
st}_{\delta^0}(\lambda)$ of $\Sigma$-stable pairs is a
compact subspace of this vanishing locus.

\begin{pr} \hfill{\break}
 i) There exists a Hausdorff open neighbourhood
$\bar{\cal B}(\Sigma)$ of ${\cal M}^{\Sigma\hbox{-}\rm
st}_{\delta_0}(\lambda)$ in $\bar{\cal B}^{**}$ such
that $Z(\hg|_{\bar{\cal B}(\Sigma)})={\cal M}^{\Sigma\hbox{-}\rm
st}_{\delta_0}(\lambda)$. In particular, ${\cal M}^{\Sigma\hbox{-}\rm
st}_{\delta^0}(\lambda)$ is an open subspace of  ${\cal
M}^{\rm ssimple}_{\delta^0}(\lambda)$.\\ 
ii) The virtual fundamental
class of ${\cal M}^{\Sigma\hbox{-}\rm
st}_{\delta^0}(\lambda)$ 
defined  by the Fredholm section $\bar\hg|_{\bar{\cal B}(\Sigma)}$
coincides with the virtual fundamental class obtained using the usual
Sobolev completion
$$\bar{\cal A}_k:=[\bar{\cal
A}(L)]^{(\lambda,v)}_{\delta^0,k}\times L^2_k(L)
$$
and the corresponding Fredhom section $\bar\hg_k$  in the Banach bundle
$$\bar\Eg_k:= [\bar{\cal A}_k^{**}]\times_{{\cal
G}^\C_{k+1}}L^2_{k-1}(\Lambda^{01}_Y\otimes L)
$$
over the Banach manifold
$$\bar{\cal B}^{**}_k:=\qmod{\bar{\cal A}^{**}_k}{{\cal G}_{k+1}^\C}\ .
$$
\end{pr}
\pf \\
$i)$ We define (compare with Definition 3.5)
$$U^s_k(\Sigma):=\{(\delta,\varphi)\in\bar{\cal A}^s_k|\
\exists\sigma\in\Sigma
\ {\rm\ such\ that}\ \varphi_j\ne 0\ \forall j\in\{1,\dots,r\}\ {\rm for\
which}\ $$
$$ \R_{\geq 0}v^j\ {\rm is\ not\ a\ face\ of}\ \sigma\}\ .
$$
The  set $U^s_k$ is obviously open in
$\bar{\cal A}^s_k$  and 
${\cal G}^\C_{k+1}$-invariant. Since $K^\C_w$ acts
freely on
$U(\Sigma)\subset\C^r$, it follows that ${\cal G}^\C_{k+1}$ acts freely
on $U^s_k(\Sigma)$, hence $U^s_k(\Sigma)\subset  [\bar{\cal A}^s_k]^{**}$.
One the other hand, an integrable pair $(\delta,\varphi)$ is obviously
$\Sigma$-stable. Put 
$$\bar{\cal B}(\Sigma):=U^s_k(\Sigma)/{\cal
G}^\C_{k+1}\ .$$
$ii) $  Since the proof uses the explicit definition of the virtual
fundamental class associated with a Fredholm section, we recall
briefly this construction. Let $s$ be a Fredholm section in a Banach
bundle $E$ over a Banach manifold $B$, and suppose, for simplicity, that
the vanishing locus $Z(s)\subset B$ is compact.

The virtual fundamental class $[Z(s)]^{\rm vir}_s$ is obtained in two
steps:\vspace{2mm}\\
a) one shows that, for a sufficiently small open neighbourhood $B_0$ of
$Z(s)$ in $B$, there exists a finite rank subbundle $E'$ of $E|_{B_0}$ such
that the  section $s_0$ induced by $s$ in the quotient bundle
$E_0:=E|_{B_0}/E'$ is regular.\vspace{2mm}\\
b) The restriction of $s$ to the finite dimensional manifold $Z(s_0)$
takes values in the bundle of finite rank $E'|_{Z(s_0)}$; let $s'$ be the
induced section in $E'|_{Z(s_0)}$. One has $Z(s)=Z(s')$, and one defines 
$$[Z(s)]^{\rm vir}_s:=[Z(s')]^{\rm vir}_{s'}\ ,
$$
where  $[Z(s')]^{\rm vir}_{s'}$ stands for the cap product
$$[Z(s')]^{\rm vir}_{s'}:=e(E'|_{Z(s_0)},s')\cap [Z(s_0)]\in
H_d(Z(s'),\Z)\ .
$$
Here $e(E'|_{Z(s_0)},s')\in H^{\rk(E')}(Z(s_0),Z(s_0)\setminus
Z(s'),\Z)$ is the localized
Euler class of $E'|_{Z(s_0)}$ with respect to the section $s'$ (see [Br], [OT2]
for details). The class  $[Z(s')]^{\rm vir}_{s'}$ is a priori an element in
the Cech homology of $Z(s')$ (which coincides with the usual
homology in our case), and it does not depend on the chosen finite rank
subbundle
$E'$.\\

Now we come back to the proof. The idea is very simple: we show  that,
applying step a) in a suitable way to the   sections $\bar\hg$,
$\bar\hg_k$ in the two bundles $\bar\Eg$, $\bar\Eg_k$,
  one gets the same    finite
dimensional manifold $Z(s_0)$, the same bundle
over $Z(s_0)$, and the same section $s'$. Therefore,   the
above two-step  procedure yields the same virtual fundamental class, because
the final step b) will be actually the same.

Consider the Banach bundle
$$\bar\Fg:=[\bar{\cal A}^s_k]^{**}\times_{{\cal G}^\C_{k+1}}
L^2_{k-1}(\Lambda^{01}_Y\otimes L)\ .
$$

One has an obvious  continuous injective  bundle morphism
$\bar\Fg\hookrightarrow
\bar\Eg$.  Note also that
$$\bar  \Fg|_{\bar{\cal B}^{**}_k}=\bar\Eg_k\ .
$$

By Corollary 5.2,    every integrable pair
$p=(\delta,\varphi)\in\bar{\cal A}^s_k$ belongs   to the
standard completed configuration space $\bar{\cal A}_k$. 
Therefore, the differential $d_p\bar h$ maps 
$$T_p(\bar{\cal A}_k)=[A^{01}(Y)^{\oplus r}]^v_k\times L^2_k(L)$$
onto $L^2_{k-1}(\Lambda^{01}_Y\otimes L)$. Let $K_p$ be the (finite
dimensional) cokernel of the restriction
$$d_p\bar h|_{[A^{01}(Y)^{\oplus r}]^v_k\times
L^2_k(L)}:[A^{01}(Y)^{\oplus r}]^v_k\times L^2_k(L)\map
L^2_{k-1}(\Lambda^{01}_Y\otimes L)\ .
$$
The family of cokernels $(K_p)_{\bar h(p)=0}$ descend to a linear space
$\Kg$ over ${\cal
M}^{\rm ssimple}_{\delta^0}(\lambda)$, which is a quotient linear space
of the restriction $\bar\Fg|_{{\cal
M}^{\rm ssimple}_{\delta^0}(\lambda)}$.

Since ${\cal M}^{\Sigma\hbox{-}\rm st}_{\delta^0}(\lambda)$ is compact,
there exists a finite system $\underline s=(s_l)_{1\leq l\leq n}$ of sections
of $\bar\Fg$ which span $\Kg_{[p]}$, for every $[p]\in  {\cal
M}^{\Sigma\hbox{-}\rm st}_{\delta^0}(\lambda)$.  A generic perturbation
$\underline s'=(s'_l)_{1\leq l\leq n}$ of $\underline s$, which is
sufficiently close to $\underline s$, will have the following properties:\\ \\
1. It still spans $\Kg_{[p]}$, for every $[p]\in  {\cal M}^{\Sigma\hbox{-}\rm
st}_{\delta^0}(\lambda)$.\\ 2. It is fibrewise linear independent in every
point   $[p]\in  {\cal M}^{\Sigma\hbox{-}\rm st}_{\delta^0}(\lambda)$.\\

The system $\underline
s'$ spans a  rank $n$ subbundle $\bar\Eg'$ of the restriction of $\bar
\Fg$ to a sufficiently small neighbourhood of ${\cal M}^{\Sigma\hbox{-}\rm
st}_{\delta^0}(\lambda)$.

But   the
inclusion
$$L^2_{k-1}(\Lambda^{01}_Y\otimes L)\map L^s_{-1}(\Lambda^{01}_Y\otimes L)
$$
induces an isomorphism  
$$\qmod{L^2_{k-1}(\Lambda^{01}_Y\otimes L)}{d_p\bar  h([A^{01}(Y)^{\oplus
r}]^v_k\times L^2_k(L))}\simeq\qmod{L^s_{-1}(\Lambda^{01}_Y\otimes
L)}{\im(d_p\bar h)}\ 
$$
for every integrable pair $p$. {This follows   using the duality
$L^s_{-1}=[L^{\frac{s}{s-1}}_1]^*$   and
 standard $L^p$ - theory.} Therefore, the system $\underline s'$ also
spans the cokernels of the absolute differentials of the section
$\bar\hg$ in the points of ${\cal M}^{\Sigma\hbox{-}\rm
st}_{\delta^0}(\lambda)$.  The section $\hg_0$ induced by
$\bar\hg$ in the quotient bundle $\bar\Eg/\bar\Eg'$ is therefore regular  on a
sufficiently small neighbourhood of  ${\cal M}^{\Sigma\hbox{-}\rm
st}_{\delta^0}(\lambda)$.  Hence one can use the finite
dimensional manifold $Z(\bar\hg_0)$ and the induced section $\bar\hg'$
in $\bar\Eg'|_{Z(\hg_0)}$ to compute the virtual fundamental class of  $[{\cal
M}^{\Sigma\hbox{-}\rm st}_{\delta^0}(\lambda)]^{\rm vir}_{\bar\hg}$.

On the other hand,   the restriction  $\bar \Eg'|_{{\cal
B}^{**}_k}$ is a  rank $n$ subbundle of $\bar \Eg_k$, and the induced
section $[\bar\hg_k]_0$ in the quotient bundle is regular around ${\cal
M}^{\Sigma\hbox{-}\rm st}_{\delta^0}(\lambda)$.

The proof is now completed   
since $\bar\Eg'\subset\bar\Fg$, so
$$Z([\bar\hg_k]_0)=Z(\bar\hg_0) 
$$
 by elliptic regularity. 
\qed

\begin{re} The method used in the proof of the  proposition above can be
used to show that  the virtual fundamental class  obtained using the standard
Sobolev completion
$\bar{\cal A}_k$, does not depend on the (sufficiently large) index $k$.

The  analogous statement should be true for all gauge
theoretical problems of  Fredholm-type.
\end{re}
\begin{pr}
  The Kobayashi-Hitchin correspondence gives an identification 
$${\cal
M}_{(\tg,g,A^0)}(\lambda)\textmap{\simeq} {\cal
M}^{\Sigma\hbox{-}\rm st}_{\bar\partial_{A^0}}(\lambda)$$
which maps $[{\cal
M}_{(\tg,g,A^0)}(\lambda)]^{\rm vir}_{(\mg_\tg,\hg)}$ onto $[{\cal
M}^{\Sigma\hbox{-}\rm st}_{\bar\partial_{A^0}}(\lambda)]^{\rm vir}_{\bar\hg}$.
\end{pr}
\pf By the Proposition above, it suffices to show that the Kobayashi-Hitchin
isomorphism maps  $[{\cal M}_{(\tg,g,A^0)}(\lambda)]^{\rm
vir}_{(\mg_\tg,\hg)}$   onto 
$[{\cal
M}^{\Sigma\hbox{-}\rm st}_{\bar\partial_{A^0}}(\lambda)]^{\rm
vir}_{\bar\hg_k}$.

But this follows by exactly same argument as in the proof of Theorem
3.2 in [OT2]. The main ingredient is the technique developed in [T] (see
also [LT]) to prove that the Kobayashi-Hitchin correspondence is a
global isomorphism of real analytic spaces.
\qed

Let $E$ be a Hermitian vector bundle on $Y$. The distribution space
$L^s_{-1} (E)$    
is the topological dual   of the Banach space
$L^t_1(E^\vee)$, where
$t$ is related to $s$ by the formula
$$\frac{1}{t}+\frac{1}{s}=1\ .
$$
Since $s<2$, it follows that $t>2$, hence on has a bounded
embedding $L^t_1( E^\vee)\subset {\cal C}^0( E^\vee)$.  In particular,
for every point $y\in Y$ one has a well defined continuous   evaluation
map
${\rm ev}_y:L^t_1(  E^\vee)\ra    E^\vee_y$.  Let
$L^t_1(  E^\vee)_y$ be the kernel of this map.

The exact sequence of Banach spaces
$$0\map L^t_1(  E^\vee)_y\map L^t_1( E^\vee)\textmap{{\rm ev}_y}
 E_y^\vee \map 0
$$
splits topologically, so one gets an exact sequence of dual spaces
$$0\map E_y \map L^s_{-1}(E )\textmap{\rho^E_y}
[L^t_1( E^\vee)_y]^{*} \map0\ .
$$

The monomorphism on
the  left is   the embedding of the space of  
Dirac distributions concentrated in $y$ in the space $L^s_{-1}(E)$ .     

More generally, if   $D=\{y_1,\dots,y_k\}\subset E$ is a finite set,
one has an exact sequence

$$0\map \bigoplus_{y\in D}E_y \map L^s_{-1}(E)\textmap{\rho^E_D}
[L^t_1( E^\vee)_D]^{*}\map  0\ ,
$$
where $L^t_1( E^\vee)_D$ is the closed subspace of $L^t_1( E^\vee)$
consisting of
$L^t_1$-sections of $ E^\vee$ which vanish   on $D$.
\begin{lm} Let $E$ be a ${\cal C}^\infty$ vector bundle on $Y$,  
$\delta$   a  semiconnection on $E$, and let ${\cal E}$ be the
corresponding holomorphic bundle.
Let $D\subset Y$  be a finite set. Consider the bounded operators
$$L^s(E)\textmap{\delta}L^s_{-1}(\Lambda^{01}_Y\otimes
E)\textmap{\rho^{[\Lambda^{01}_Y\otimes
E]}_D} [L^t_1(\Lambda^{01}_Y\otimes  E)^\vee_D]^{*}\ :
$$
\begin{enumerate}
\item  One has
$$\ker \left[\rho^{[\Lambda^{01}_Y\otimes
E]}_D\circ \delta\right]=H^0({\cal E}(D))\ .$$
\item When $\deg( E^\vee\otimes K_Y(-D))<0$, the map
$\rho^{[\Lambda^{01}_Y\otimes E]}_D\circ \delta$ is surjective.  
\end{enumerate}
\end{lm} 
\pf\\
1. A section $\varphi\in L^s(E)$ belongs to $\ker
\left[\rho^{[\Lambda^{01}_Y\otimes E]}_D\circ \delta\right]$ if and only
if $\delta\varphi$ is a linear combination of Dirac distributions
concentrated in the points of $D$.  This means that $\varphi$ is
meromorphic with poles of order at most 1 in the points of $D$.
\\
2. The composition $\rho^{[\Lambda^{01}_Y\otimes E]}_D\circ \delta$ is the
adjoint of the composition
$$L^t_1([\Lambda^{01}_Y]^\vee\otimes  
E^\vee)_D\stackrel{i_D}\hookrightarrow
L^t_1([\Lambda^{01}_Y]^\vee\otimes   E^\vee) 
\textmap{\delta'}L^t_1(  E^\vee)\ ,
$$
where the operator $\delta'$ is the adjoint of $\delta$, i. e. it
satisfies the identity
$$\int_Y\langle \delta' (\alpha),\beta\rangle vol_g=\int_Y\langle  
\alpha,\delta(\beta)\rangle vol_g
$$
for all smooth sections $\alpha\in A^0([\Lambda^{01}_Y]^\vee\otimes  
E^\vee)$, $\beta\in A^0( E)$.

The crucial observation is that $\delta'$ is just  the tensor product
semiconnection $\delta_{\rm can}\otimes \delta^\vee$ on the bundle
$\Lambda^{10}_Y\otimes E^\vee$, where
$\delta_{\rm can}$ is the canonical holomorphic structure on
$\Lambda^{10}_Y=K_Y$ and $\delta^\vee$ is the dual semiconnection on $E^\vee$.

On the other hand, the composition $\delta'\circ i_D$ is Fredholm,
because the space  $L^t_1([\Lambda^{01}_Y]^\vee\otimes  
E^\vee)_D$ is closed and has finite codimension in
$L^t_1([\Lambda^{01}_Y]^\vee\otimes   E^\vee)$.

Therefore the surjectivity of the composition
$\rho^{[\Lambda^{01}_Y\otimes E]}_D\circ \delta$ is equivalent to the
injectivity of $\delta'\circ i_D$. Let ${\cal E}$ be the holomorphic
bundle defined by $\delta$.

By elliptic regularity we have  
$$\ker(\delta'\circ i_D)=\{s\in H^0(K_Y\otimes {\cal E}^\vee)|\
s|_D=0\}=H^0(K_Y\otimes {\cal E}^\vee(-D))\ ,
$$
which vanishes, when $\deg( E^\vee\otimes K_Y(-D))<0$.
\qed

Now we come back to our gauge theoretical problem. 
For  $1\leq j\leq r$ let
$D_j \subset Y$ be a finite set of points,  and denote by $D$ the system 
$(D_j)_j$. The same argument applied to the line bundles
$\Lambda^{01}_Y\otimes L_j$ yields the exact sequence
$$0\map \bigoplus_j\bigoplus_{y\in D_j}(\Lambda^{01}_Y\otimes
L_j)_{y}\map
\bigoplus_jL^s_{-1}(\Lambda^{01}_Y\otimes L_j)\textmap{\rho_D}
\bigoplus_j [L^t_1(\Lambda^{01}_Y\otimes L_j)_{D_j}^\vee]^{*}\map 0\ ,
$$
where $L^t_1(\Lambda^{01}_Y\otimes L_j)_{D_j}$ is the space of $L^t_1$-sections of
$\Lambda^{01}_Y\otimes L_j$ which vanish on $D_j$.

We denote by $\bar\Eg_0$  the associated bundle of the
principal bundle $[\bar{\cal A}^s_k]^{**}\ra \bar{\cal B}^{**}$ with
standard fibre 
$$L^s_{-1}(\Lambda^{01}_Y\otimes L)_D:=\bigoplus_j
[L^t_1(\Lambda^{01}_Y\otimes L_j)_{D_j}^\vee]^{*}\ ,$$
and by $\bar\Eg'$ the   associated bundle with fibre 
$$ (\Lambda^{01}_Y\otimes
L)_{D}:=\bigoplus_j\bigoplus_{y\in D_j}(\Lambda^{01}_Y\otimes
L_j)_{y} \ .
$$
One obtains an exact sequence of Banach bundles over ${\cal B}^{**}$:
$$0\map\bar \Eg'\map\bar \Eg\stackrel{\rg_D}{\map}\bar\Eg_0\map 0\ ,
\eqno{(*)}
$$
where ${\rg_D}$ is the bundle epimorphism induced by $\rho_D$. Let
 $\bar\hg_0$ be the induced section $\bar\hg_0:=\rg_D\circ\bar\hg$.

\begin{pr}  Suppose that $\#D_j+\deg(L_j)> 2g_Y-2$ for all
$j\in\{1,\dots,r\}$. Then   the section $\bar\hg_0:=\rg_D\circ\bar\hg$ is
regular at every point
$[\delta,\varphi]$ of its vanishing locus.
\end{pr}
\pf It   suffices to show that the differential at
$p=(\delta,\varphi)$ of the restriction of
$\rho_D\circ \bar h$ to the slice $V_{p,\varepsilon_p}$ provided by
Proposition 5.3. is surjective.\footnote{In the infinite dimensional
framework the surjectivity of the differential at a point does not
suffices to assure that the map is a submersion at that point. One also
has to check that the kernel of the differential has a closed
complement.  But this condition is obviously satisfied in our case, since
this kernel has finite dimension.}

On the other hand, the differential of $\bar h$ at $p$ vanishes on the
tangent space at the orbit $p\cdot{\cal G}_{k+1}^\C$, so it suffices to
show that the differential of $\rho_D\circ \bar h$ at $p$ is surjective.
But
$$\left.\frac{\partial}{\partial \varphi}\right|_p(\rho_D\circ \bar
h)=\rho_D\circ\delta\ ,
$$
so the result follows from the previous lemma.
\qed

We recall   Brussee's associativity principle for   virtual
fundamental classes associated with Fredholm sections ([Br], [OT2]) in the
special case of sections with compact vanishing locus.
\begin{thry}
Let
$$0\map E'\map E\map E_0\map 0
$$
be an exact sequence of Banach bundles over a Banach manifold $B$, and let
$s$ be a Fredholm section in $E$ with compact vanishing locus $Z(s)$.
Suppose that the induced section $s_0$ in the bundle $E_0$ is regular in
every point of its vanishing locus $Z(s_0)$. The restriction $s|_{Z(s_0)}$ can be
regarded as a section
$s'$ in $E'|_{Z(s_0)}$. Via the obvious identification $Z(s)=Z(s')$ one has
$$[Z(s)]^{\rm vir}_s=[Z(s')]^{\rm vir}_{s'}\ .
$$

\end{thry}

Using this we   can prove our main result  
\begin{thry} The identifications given by the Kobayashi-Hitchin
correspondence (Theorem 3.6) and by the embedding theorem (Theorem 4.6)
map  the virtual fundamental class induced by the Fredholm
description of the moduli space ${\cal M}_{(\tg,g,A^0)}(\lambda)$ onto
the algebraic geometric virtual fundamental class   defined by the system of
sections
$e^i_j$ on the smooth algebraic variety
$X_{\Sigma,P'}( p_{\underline c'}^*({\Lg}^{{\underline
c'}}_{y_0})|_{P'})$.
\end{thry}
\pf First note that, by Proposition  5.6, it  suffices to compare the
gauge theoretical virtual fundamental class $[{\cal
M}^{\Sigma\hbox{-}\rm st}_{\delta^0}(\lambda)]^{\rm vir}_{\bar\hg}$
with the algebraic geometric virtual fundamental class defined by the
sections $e^i_j$ on the smooth algebraic variety $X_{\Sigma,P'}(
p_{\underline c'}^*({\Lg}^{{\underline c'}}_{y_0})|_{P'})$.

We apply the associativity principle to the restriction of the exact
sequence
$(*)$ to the open  set $\bar{\cal B}(\Sigma)$ of $\bar{\cal B}^{**}$. 
The hypothesis of this principle is verified by  Proposition 5.8. 

Put $Z_0:=Z(\bar\hg_0|_{\bar{\cal B}(\Sigma)})$, and let
$\bar\hg'$ be the induced section in $\Eg'|_{Z_0}$. It follows that 
$$[{\cal
M}^{\Sigma\hbox{-}\rm st}_{\delta^0}(\lambda)]^{\rm
vir}_{\bar\hg}=[Z(\bar\hg')]^{\rm
vir}_{\bar\hg'}\ . 
$$
We claim that:\vspace{2mm}\\ 
a) With the notations of section 4, there is a natural identification
$$Z_0=[\iota_{\underline D}|_P]^*(X_{\Sigma,P'}(
p_{\underline c'}^*({\Lg}^{{\underline c'}}_{y_0})|_{P'}))\ ,
$$
where, on the right $X_{\Sigma,P'}(
p_{\underline c'}^*({\Lg}^{{\underline c'}}_{y_0})|_{P'})$ was
considered as a toric fibre bundle over $P':=\iota_{\underline D}(P)$.
\\
b) Via the identification above, the sections $[\iota_{\underline
D}|_P]^*(e^i_j)$ coincide with the components of $\bar\hg'$.\\

The proof of the two claims is   easy: for a)  one uses the
first statement in Lemma 5.7 to get a set theoretical bijection; then
the method used in the proof of Theorem 4.5 gives the needed
isomorphism. Note that this time, one just has to identify two smooth
complex manifolds, so it is not necessary to take   
ringed space  structures into account. For b) it suffices to notice that the
components of $\bar\hg'$ and the sections $[\iota_{\underline
D}|_P]^*(e^i_j)$ are induced by evaluation maps associated with the points of
$D_j$.

On the other hand,  by the Lemma below, 
Brussee's virtual fundamental class associated with a regular section in
an algebraic vector bundle over a smooth algebraic variety coincides with
the corresponding Fulton virtual fundamental class via the cycle map (see
[Fu], ch. 19.1 - 19.2).
\qed
\begin{lm} Let $E\ra X$ be an algebraic vector bundle of rank $r$ over a
smooth
$n$-dimensional algebraic variety $X$, and let $s$ be a section in $E$.

Then the Fulton virtual fundamental class $\Z(s)\in A_{n-r}(Z(s))$ is mapped
onto the Brussee virtual fundamental class $[Z(s)]^{\rm vir}_s$ via the
cycle map
$$cl_{Z(s)}:A_{n-r}(Z(s))\map H_{n-r}^{\rm BM}(Z(s),\Z)\ ,
$$
where $A_{*}(Z(s))$ denotes the  Chow groups  of the complex scheme $Z(s)$ and
$H_{*}^{\rm BM}(Z(s),\Z)$   the Borel-Moore homology groups of $Z(s)$.
\end{lm}
\pf   We will give the proof in
the case when
$Z(s)$ is compact. This case is sufficient for our purposes, and the general
case follows the same idea. 

Note first that the Cech homology of $Z(s)$, to which $[Z(s)]^{\rm
vir}_s$ belongs, coincides with the standard
 singular homology ([Br], Ch. 2). Let $s_E:X\ra E$  be the
zero-section of $E$. By definition, one has
$$[Z(s)]^{\rm vir}_s=s^*(\tau_E)\cap [X]\in H_{n-r}(Z(s),\Z)\ ,
$$
where $\tau_E\in H^r(E,E\setminus\im(s_E),\Z)$ is the Thom class of $E$, and
$[X]$ is the fundamental class of $X$ in Borel-Moore homology.

 Since $E$ is a smooth
algebraic variety, the class
$\Z(s)$ can be identified with the {\it refined intersection  product} of
algebraic cycles
$\im(s_E)\cdot\im(s)\in A_{n-r}(\im(s_E)\cap\im(s))$ via the  obvious
identification $\im(s_E)\cap\im(s)=Z(s)$ (see [Fu] ch. 8.1, Corollary 8.1.1,
Ch. 14.1).  Therefore, by   Corollary 19.2 [Fu], we get
$$cl_{Z(s)}(\Z(s))= cl_{Z(s)}(\im(s_E)\cdot\im(s))=
cl_{Z(s)}(\im(s_E))\cdot
cl_{Z(s)}(\im(s))=$$
$$=[\im(s_E)]\cdot[\im(s)]\ ,
$$
where the  dot   in the last two  term stands for the refined topological
intersection ([Fu], p. 378).  Note now that the Poincar\'e dual
$$PD([\im(s_E)])\in H^r(E,E\setminus \im(s_E),\Z)$$
of $[\im(s_E)]$ in $E$ is just the Thom class $\tau_E$.  Therefore, 
$$[\im(s_E)]\cdot [\im(s)]=\tau_E\cap s_*[X]  =
(s|_{Z(s)})_*\left( s^*(\tau_E)\cap[X]\right) 
$$
in $H_*(\im(s_E)\cap \im(s),\Z)$. This shows that
$$cl_{Z(s)}(\Z(s))=  s^*(\tau_E)\cap[X]=[Z(s)]^{\rm vir}_s\ .
$$
\qed

 {\small
Authors addresses: \vspace{2mm}\\
Ch. Okonek,   Institut f\"ur Mathematik, Universit\"at Z\"urich,  
Winterthurerstrasse 190, CH-8057 Z\"urich, Switzerland,  e-mail:
okonek@math.unizh.ch\vspace{1mm}\\
 A. Teleman, LATP, CMI,   Universit\'e de Provence,  39  Rue F. J.
Curie, 13453 Marseille Cedex 13, France,  e-mail: teleman@cmi.univ-mrs.fr 
    ,  and\\ Faculty of Mathematics, University of Bucharest, Bucharest,
Romania }


\begin{thebibliography}{CGMS}

\bibitem[A]{A} Audin, M.: {\it The topology of torus actions on symplectic
manifolds}, Progress in Math. 93, Birkh\"auser, Boston-Basel-Berlin, 1991.

\bibitem[B]{B} Batyrev, V.: {\it Quantum cohomology rings of toric manifolds},
Journ\'ees de G\'eomt\'erie Alg\'ebrique d'Orsay (Juillet 1992), Ast\'erisque
218, 9-34, Soci\'et\'e Math\'ematique de France, Paris 1993.

\bibitem[BF]{BF} Behrend, K., Fantechi, B.: {\it The intrinsic normal cone},
Invent. Math. 128, 1,    45-88, 1997.
 

\bibitem[Br]{Br}  Brussee, R.: {\it The canonical class and the ${\cal C}^\infty$
properties of K\"ahler surfaces}, New York J. Math. 2 , 103-146, 1996.

\bibitem[Bra]{Bra} Bradlow, S., B.: {\it Special metrics and
stability for holomorphic bundles with global
sections},   J. Diff. Geom., 33, 169-214,  1991.
 

\bibitem[CGMS]{CGMS}   Cieliebak, K., Gaio,  A. R., Salamon,  D., Mundet i Riera,
I.: {\it  The symplectic vortex equations and invariants of Hamiltonian group
actions}, {\tt math.SG/0111176}.

\bibitem[CGS]{CGS}   Cieliebak, K., Gaio,  A. R., Salamon,  D.:
{\it J-holomorphic curves, moment maps, and invariants of Hamiltonian group
actions}, Internat. Math. Res. Notices no 16, 831-882, 2000. 

\bibitem[CK]{CK} Cox, D., Katz, S.: {\it Mirror symmetry and algebraic
geometry}, Mathematical Surveys and Monographs 68, AMS Providence, RI, 1999.

\bibitem[CKL]{CKL} Cox, D., Katz, S., Li, Y.: {\it Virtual fundamental classes
of zero loci}, Contemporary Mathematics 276,  157-166, 2001.


\bibitem[Co]{Co} Cox, D.: {\it Recent developments in toric geometry},
Proceedings of Symposia in Pure Mathematics Vol. 62, 2, 389-436, 1997.


\bibitem[D1]{D1}  Donaldson, S. K.: {\it Anti self-dual
Yang-Mills connections over complex algebraic
surfaces and stable vector bundles}, Proc.
London Math. Soc., 3,  50,  1-26, 1985.

\bibitem[D2]{D2}  Donaldson, S. K.: {\it Irrationality  and 
the\ $h$-cobordism\ conjecture}, J.
Diff. Geom., 26, 141-168, 1987.


\bibitem[DK]{DK}  Donaldson, S.,   Kronheimer, P.: {\it The Geometry of
Four-Manifolds}, Oxford Univ. Press, 1990.

\bibitem[Dou]{Dou}  Douady, A.: {\it Le probl\`eme des modules
pour les vari\'et\'es analytiques complexes.},  Sem.
Bourbaki   nr. 277, 17-eme annee 1964/65.

\bibitem[DOT]{DOT} D\"urr, M., Okonek, Ch., Teleman, A.: {\it From
Poincar\'e to Seiberg-Witten},   in preparation.



\bibitem[D\"u]{Du}  D\"urr, M.: {\it Virtual fundamental classes and Poincar\'e
invariants},    in preparation.

\bibitem[F]{F} Fulton, W.:{\it Intersection Theory}, Springer Verlag,
New York-Berlin-Heidelberg, 1984.

\bibitem[G]{G}   Gaio, A. R.  : {\it J-holomorphic curves and moment maps},  
Ph. D. Thesis, University of Warwick, November 1999.

\bibitem[Gi1]{Gi1} Givental, A.: {\it Equivariant Gromov-Witten invariants},
Internat. Math. Res. Notices, 13, 613-663, 1996.

\bibitem[Gi2]{Gi2} Givental, A.: {\it A symplectic fixed point theorem for
toric manifolds}, The Floer memorial volume,  Progr. Math. 133,
Birkh\"auser, Basel, 445-481, 1995.
\bibitem[GS]{GS}   Gaio, A. R., Salamon,  D.: {\it Gromov-Witten  invariants of
symplectic quotients and adiabatic limits}, {\tt math.SG/0106157}. 


 

 

\bibitem[Gr]{Gr}   Gromov, M.: {\it Pseudo-holomorphic curves in symplectic
manifolds},   Invent. math.   82,   307-347,  1985.

\bibitem[GP]{GP} Graber, T.,  Pandharipande, R.: {\it   Localization of
virtual  classes}, Invent. math. 135, 487--518, 1999. 

\bibitem[Ha1]{Ha1}  Halic, M.: {\it Higher genus curves on toric
varieties}, preprint, {\tt math.AG/0101100}. 

\bibitem[Ha2]{Ha2} Halic, M.: {\it Hamiltonian Gromov-Witten invariants
of toric varieties}, in preparation.

\bibitem[Hi]{Hi} N.\ Hitchin, {\it The self-duality equations on a Riemann 
surface}\rm ,  Proc.  London  Math. Soc. , 55, 59-126,  1987.

\bibitem[K]{K} Kirwan, F.: {\it Cohomology of quotients in symplectic and
algebraic geometry}, Princeton Univ. Press, Princeton, NJ, 1984.


\bibitem[LiT]{LiT} Li, J., Tian, G.: {\it  Virtual moduli cycles and
Gromov-Witten invariants of general symplectic manifolds}, Topics in symplectic
$4$-manifolds (Irvine, CA, 1996), 47 - 83, First Int. Press Lect. Ser., I,
Internat. Press, Cambridge, MA, 1998.

\bibitem[LL]{LL} L\"ubke, M.,  Lupa\c{s}cu, P.: {\it Isomorphy of the
gauge theoretical and the deformation theoretical moduli space of
simple holomorphic pairs},   preprint, Leiden, 2001.

\bibitem[Lu]{Lu}  Lupa\c{s}cu, P.: {\it Seiberg-Witten equations and complex
surfaces},   Ph. D. Thesis, Z\"urich, 1998.

\bibitem[LO]{LO} L\"ubke, M.,  Okonek, Ch.: {\it Moduli
spaces of simple bundles and
Hermitian-Einstein  connections},    Math.
Ann.   276, 663-674,  1987.

\bibitem[LOT]{LOT} L\"ubke, M., Okonek, Ch., Teleman, A.: {\it Moduli
spaces in gauge theory and algebraic geometry},   in preparation.

\bibitem[LT1]{LT1} L\"ubke, M., Teleman, A.: {\it The Kobayashi-Hitchin 
correspondence},    World Scientific Publishing Co.,  1995.

\bibitem[LT2]{LT2} L\"ubke, M., Teleman, A.: {\it The universal
Kobayashi-Hitchin  correspondence for framed pairs},    in preparation.


\bibitem[MP]{MP} Morrison, D., Plesser, M.:{\it Summing the instantons: quantum
cohomology and mirror symmetry in toric varieties}, Nuclear Phys. B, 440, no
1-2, 279-354, 1995.

\bibitem[Mu1]{Mu1}  Mundet i Riera, I.: {\it
A Hitchin--Kobayashi correspondence for Kaehler fibrations,}
 {J. Reine Angew. Math.}, 528, 41-80, 2000.



\bibitem[Mu2]{Mu2}   Mundet i Riera, I.: {\it Hamiltonian Gromov--Witten 
invariants}, {\tt math. SG/0002121}.

\bibitem[OP]{OP} Oda, T., Park, H.: {\it Linear Gale transforms and
Gelfand-Kapranov-Zelevinsky decompositions}, T\^ohoku Math. J. 43, 375-395,
1991.

\bibitem[OT1]{OT1}    Okonek, Ch.,  Teleman, A.: {\it Quaternionic
monopoles}, Commun.\ Math.\ Phys., Vol. 180, Nr. 2, 363-388, 1996. 
 

\bibitem [OT2]{OT2}   Okonek, Ch., Teleman, A.: {\it Gauge theoretical
equivariant Gromov-Witten invariants  and the full Seiberg-Witten
invariants of ruled surfaces},    to appear in Commun. Math. Phys.

 
 
\bibitem [R]{R} Ruan, Y.: {\it Topological sigma model and Donaldson-type
invariants in Gromov theory},  Duke Math. J., Vol. 83,   no. 2, 461- 500,
1996.

\bibitem [S]{S} Serre, J. P.: {\it  G\'eometrie alg\'ebrique et G\'eometrie
analytique}, Ann. Inst. Fourier 6, 1-42,  1956.

\bibitem [Sp]{Sp} Spielberg, H.: {\it  The Gromov-Witten invariants of
symplectic toric manifolds},  {\tt mathAG/0006156}.

 
\bibitem[Su]{Su} Suyama, Y.: {\it The analytic moduli space of
simple  framed holomorphic pairs}, Kyushu J. Math., Vol. 50,   65-68, 
1996.

\bibitem[Te]{Te} Teleman, A.: {\it Stable vector bundles over
non-K\"ahlerian surfaces}, Ph. D. Thesis, Z\"urich, 1993.


\bibitem[Th1]{Th1} Thaddeus, M.: {\it Stable pairs, linear  systems and the Verlinde
formula}, Invent. math. 117, 181-205,  1994.


\bibitem[Th2]{Th2}   Thaddeus, M.: {\it Geometric invariant theory
 and flips}, 
J. AMS, Vol. 9 , 691-725, 1996.

 


\bibitem[W1]{W1}  Witten, E.: {\it The Verlinde algebra and the cohomology of the
Grassmannian},  Geometry, Topology  and  Physics, 357 - 422, Conf. Proc. Lecture
Notes Geom. Topology, IV, Internat. Press, Cambridge, MA, 1995.

\bibitem[W2]{W2} Witten, E.: {\it Monopoles and four-manifolds}, Math.  Res.
Letters 1,  769-796, 1994.

\bibitem[W3]{W3} Witten, E.: {\it Phases of $N=2$ theories in two dimensions},
Nuclear Physics B, 403, no 1-2, 159-222, 1993.

\bibitem[W4]{W4} Witten, E.: {\it Two-dimensional gravity and intersection
theory on moduli spaces}, Survey in Diff. Geom. 1, 243-310, 1991.



 \end{thebibliography}
\end{document}